\documentclass[11pt]{article}
\usepackage{color}

\usepackage{amsmath,amssymb,mleftright,mathtools}
\usepackage{graphicx}
\usepackage{theorem}
\usepackage{here}

\newtheorem{thm}{Theorem}[section]
\newtheorem{lem}{Lemma}[section]
\newtheorem{cor}{Corollary}[section]
\newtheorem{prp}{Proposition}[section]
\theorembodyfont{\rmfamily}

\makeatletter

\@addtoreset{equation}{section}
\makeatother

\def\L{{{\cal L}}}

\def\al{{\alpha}}

\def\de{{\delta}}

\def\la{{\lambda}}
\def\si{{\sigma}}

\def\th{{\theta}}

\def\bde{{\text{\boldmath $\delta$}}}

\def\bth{{\text{\boldmath $\theta$}}}

\def\bxi{{\text{\boldmath $\xi$}}}

\def\bmu{{\text{\boldmath $\mu$}}}

\def\bnabla{{\text{\boldmath $\nabla$}}}
\def\bel{{\text{\boldmath $\ell$}}}
\def\alh{{\widehat \al}}

\def\Th{{\Theta}}
\def\De{{\Delta}}
\def\Si{{\Sigma}}

\def\La{{\Lambda}}

\def\Deh{{\widehat \De}}
\def\bTh{{\text{\boldmath $\Th$}}}

\def\bSi{{\text{\boldmath $\Si$}}}

\def\bLa{{\text{\boldmath $\La$}}}

\def\bPsi{{\text{\boldmath $\Psi$}}}

\def\bThh{{\widehat \bTh}}
\def\bSih{{\widehat \bSi}}
\def\bPsih{{\widehat \bPsi}}

\def\x{{\text{\boldmath $x$}}}

\def\A{{\text{\boldmath $A$}}}
\def\B{{\text{\boldmath $B$}}}
\def\C{{\text{\boldmath $C$}}}

\def\G{{\text{\boldmath $G$}}}
\def\H{{\text{\boldmath $H$}}}
\def\I{{\text{\boldmath $I$}}}

\def\O{{\text{\boldmath $O$}}}
\def\P{{\text{\boldmath $P$}}}

\def\R{{\text{\boldmath $R$}}}
\def\S{{\text{\boldmath $S$}}}

\def\U{{\text{\boldmath $U$}}}
\def\V{{\text{\boldmath $V$}}}
\def\W{{\text{\boldmath $W$}}}
\def\X{{\text{\boldmath $X$}}}
\def\Y{{\text{\boldmath $Y$}}}

\def\xb{{\overline \x}}

\def\Xb{{\overline \X}}

\def\ah{{\hat a}}
\def\bh{{\hat b}}

\def\Rh{{\widehat R}}

\def\Dc{{\cal D}}

\def\Nc{{\cal N}}

\def\Re{{\mathbb{R}}}
\def\tr{{\rm tr\,}}
\def\diag{{\rm diag\,}}

\def\[{{\text{\boldmath $[$}}}
\def\]{{\text{\boldmath $]$}}}

\def\zero{{\bf\text{\boldmath $0$}}}

\def\|{{\,|\,}}

\def\/{{\Bigr/\!\!}}

\def\1r{{\rm (1)}}
\def\2r{{\rm (2)}}
\def\3r{{\rm (3)}}
\def\4r{{\rm (4)}}
\def\5r{{\rm (5)}}

\def\non{{\nonumber}}
\def\Bbh{{\widehat \B}}

\renewcommand{\Im}{\mathrm{Im}}
\newcommand{\as}{\text{\ a.s.}}
\newcommand{\asl}{\text{\ as\ }}

\newcommand{\trRS}[1]{\tr[(\W+\alh\I)^{-#1}\W]}
\newcommand{\Rid}[1]{(\W+\alh\I)^{-#1}}

\def\L{{\text{\boldmath $L$}}}

\begin{document}
\title{Ridge-type Linear Shrinkage Estimation of the Matrix Mean of High-dimensional Normal Distribution}
\author{Ryota Yuasa\footnote{Graduate School of Economics, University of Tokyo, 
7-3-1, Hongo, Bunkyo-ku, Tokyo 113-0033, Japan \quad
E-Mail: ryooyry@yahoo.co.jp}
\quad and \quad
Tatsuya Kubokawa\footnote{Faculty of Economics, University of Tokyo, 
7-3-1, Hongo, Bunkyo-ku, Tokyo 113-0033, Japan \quad
E-Mail: tatsuya@e.u-tokyo.ac.jp }}
\maketitle
\begin{abstract}
The estimation of the mean matrix of the multivariate normal distribution is addressed in the high dimensional setting.
Efron-Morris-type linear shrinkage estimators based on ridge estimators for the precision matrix instead of the Moore-Penrose generalized inverse are considered, and the weights in the ridge-type linear shrinkage estimators are estimated in terms of minimizing the Stein unbiased risk estimators under the quadratic loss. 
It is shown that the ridge-type linear shrinkage estimators with the estimated weights are minimax, and that the estimated weights converge to the optimal weights in the Bayesian model with high dimension by using the random matrix theory.
The performance of the ridge-type linear shrinkage estimators is numerically compared with the existing estimators including the Efron-Morris and James-Stein estimators.

\par\vspace{4mm}
\end{abstract}

\section{Introduction}
\label{sec:int}

Statistical inference with high dimension has received much attention in recent years, because high-dimensional data have been handled in many research areas including genomics and financial engineering.
Especially, in multivariate analysis with a dimension larger than a sample size, the sample covariance matrix is ill-conditioned, and we cannot employ the inverse of the sample covariance matrix for estimating the precision matrix, which appears in multivariate statistical inference such as classification, testing hypothesis, confidence region, regression analysis and others.
For the estimation of the large covariance matrix, see Fan, Liao and Mincheva (2013), Ledoit and Wolf (2012, 15), Bodnar, Gupta and Parolya (2014), Wang, Pan, Tong and Zhu (2015), Ikeda and Kubokawa (2016) for example.
A similar difficulty appears in estimation of a mean matrix of a multivariate normal distribution.
The shrinkage estimators improving on the sample mean of the mean matrix have been provided by Efron and Morris (1972, 76), Konno (1991, 92), Tsukuma and Kubokawa (2007, 15) as an extension of the Stein phenomenon in estimation of a normal mean vector.
Using the Bayesian argument, Efron and Morris (1976) demonstrated that the estimation of the mean matrix can be reduced to the estimation of a precision matrix.
Recently, Tsukuma and Kubokawa (2015) gave a unified theory of minimaxity in the estimation of the mean matrix in high and low dimensions.
However, their estimators were limited to shrinkage estimators based on the Moore-Penrose generalized inverse, and their attention focused on the minimaxity of the shrinkage estimators.
In this paper, we consider classes of linear shrinkage estimators of the mean matrix using ridge estimators for the precision matrix instead of the Moore-Penrose generalized inverse.
We derive estimators of the optimal weights in light of minimizing the Stein unbiased risk estimators (SURE), and provide conditions for minimaxity of the linear shrinkage estimators with the estimated weights.
We also demonstrate that the estimated weights converge to the optimal weights using the random matrix theory in high dimension.

\medskip
Let us consider the model that the $p$-variate random vectors $\x_1,\ldots,\x_n$ are mutually independent and each $\x_i$ is distributed as the multivariate normal distribution
\begin{equation}
\x_i \mid \bth_i \sim \Nc_p(\bth_i, \bSi) \hspace{3mm}(i=1,\cdots, n),\label{defx}
\end{equation}
where $\bSi$ is a known and positive definite covariance matrix.
Let $\X=(\x_1\ldots,\x_n)$ and $\bTh=(\bth_1\ldots,\bth_n)$.
Consider the problem of estimating the $p\times n$ mean matrix $\bTh$ by estimator $\bThh$ relative to the quadratic loss function
\begin{equation}
L(\bTh,\bThh)={1\over np}\tr[(\bThh-\bTh)^\top  \bSi^{-1}(\bThh-\bTh)],\label{loss}
\end{equation}
When $n>p+2$, as empirical Bayes estimators, Efron and Morris (1972, 1976) proposed the shrinkage estimators
\begin{align}
\bThh^{EM1}=&\X - (n-p-2)\W^{-1}(\X-\Xb)
, \label{em1}
\\
\bThh^{EM2}=&\X - \Big\{(n-p-2)\W^{-1}+{p^2+p-2\over \tr[\W]}\I\Big\}(\X-\Xb)
,\label{em2}
\end{align}
and showed that $\bThh^{EM1}$ is minimax and improved on by $\bThh^{EM2}$ for $p>2$, where 
\begin{align}
\W&=(\X-\Xb)(\X-\Xb)^\top  \bSi^{-1}\quad {\rm and}\quad \Xb=\xb\text{\bf 1}_n^\top, \label{W}
\end{align}
for $\xb=n^{-1}\sum_{j=1}^n\x_j$ and $\text{\bf 1}_n^\top  =(1,\ldots,1)$.
However, these estimators are not available when the dimension $p$ is greater than $n$, because $\W$ is singular. Also these are likely to be unstable when the dimension $p$ is close to $n$, because $\W$ is ill-conditioned.
An alternative method is the generalized inverse $\W^-$, and Tsukuma and Kubokawa (2015) suggested Efron-Morris-type estimators using the Moore-Penrose generalized inverse $\W^+$.

\medskip
In this paper, we use the ridge-type inverse instead of the generalized inverse and consider the linear shrinkage estimators
\begin{align}
\bThh_{a,0}^{RLS}=& \X -a\Rid{1}(\X-\Xb)
, \label{R1}
\\
\bThh_{a,b}^{RLS}=&\X -\big\{a\Rid{1}+b(\tr\W)^{-1}\I\big\}(\X-\Xb)
, \label{R2}
\end{align}
where $a$ and $b$ are scalars and $\alh$ is a ridge function of $\W$.
We call these estimators {\it Ridge-type Linear Shrinkage Estimators}.
It is noted that the ridge-type inverse $\Rid{1}$ exists as long as $\alh$ is not zero.
For $\alh$, we here treat the two cases that $\alh=c$ and $\alh=c\tr[\W]$ for positive constant $c$.

\medskip
The scalar weights $a$ and $b$ control the shrinkage weight between $\X$ and $\Xb$, and it is important how to choose the weights.
In the finite sample setting, we suggest to estimate the shrinkage weights in terms of minimizing the Stein unbiased risk estimators (SURE).

\medskip
When the shrinkage weights $a$ and $b$ are functions of $\X$, it is harder to establish minimaxity of the ridge-type linear shrinkage estimators, and conditions for the minimaxity have not been provided.
One of the major contributions in this paper is to derive sufficient conditions for the minimaxity of the ridge-type linear shrinkage estimators with the estimated weights, where the ridge function is $\alh=c$ or $\alh=c\tr[\W]$.

\medskip
The other contribution is to show that the estimated weights of $a$ and $b$ converge to the optimal weights.
In estimation of a high-dimensional mean vector, the derivation of the optimal weights in linear shrinkage estimators has been studied by Wang et al.(2014), Xie et al.(2016), Bodnar et al.(2019) and others.
Using a similar argument, in the setting of the Bayesian model, we derive the optimal $a$ and $b$ in terms of minimizing the risk function of the estimators $\bThh_{a,0}^{RLS}$ and $\bThh_{a,b}^{RLS}$ under the quadratic loss.
Then, using the method of the random matrix theory, we show that the estimated weights converge to the optimal weights.

\medskip
The paper is organized as follows:
In Section \ref{sec:stein}, we treat the estimation of the mean matrix from the frequentist's point of view. We introduce some features of ridge-type linear shrinkage estimator. We obtain several sufficient conditions for the minimaxity of the ridge-type linear shrinkage estimators with known weights by the SURE method.
In Section \ref{sec:PWMstein}, we also derive the estimators of the weights $a$ and $b$ by the SURE method. We also obtain several sufficient conditions for the minimaxity of the ridge-type linear shrinkage estimators with the estimated weights.
In Section \ref{sec:RMT}, using the argument as in Efron and Morris (1972), we consider the Bayesian model and reduce the estimation of the mean matrix to the problem of estimating a precision matrix.
In the high-dimensional setting, we show that the estimated weights converge to the optimal weights of $a$ and $b$ based on the random matrix theory.
In Section \ref{sec:sim}, the performances of the ridge-type linear shrinkage estimators are numerically compared with the existing estimators including the Efron-Morris and James-Stein estimators and the low-rank estimator suggested by Gavish and Donoho (2017), and it is shown that the proposed estimators have good performances.
All the proofs are given in the Appendix.

\section{Proposal Weights in Ridge-type Linear Shrinkage Estimators and Their Minimaxity}
\label{sec:stein}

\subsection{Ridge-type linear shrinkage estimation}

Consider the linear shrinkage estimator
\begin{equation}\label{LS0}
\bThh^{LS}(\B) = \X - \B(\X-\Xb),
\end{equation}
where $\B$ is a $p\times p$ unknown matrix such that $\bSi^{-1/2}\B\bSi^{1/2}$ and $\I_p-\bSi^{-1/2}\B\bSi^{1/2}$ are symmetric and nonnegative definite.
This estimator is motivated from the Bayes estimator for the prior $\bth_i \mid \bmu \sim \Nc_p(\bmu, \bLa)$ and $\bmu\sim {\rm Uniform\  on\ } \Re^p$.
In fact, the Bayes estimator can be written as $\bThh^{Bays}=\X-\bSi(\bSi+\bLa)^{-1}(\X-\Xb)$, which yields (\ref{LS0}) for $\B=\bSi(\bSi+\bLa)^{-1}$.
The risk function of $\bThh^{LS}(\B)$ relative to the quadratic loss (\ref{loss}) is
\begin{align*}
E[L(\bTh,\bThh^{LS})]=&{1\over np}E[ \tr\{\bSi^{-1}\B(\X-\Xb)(\X-\Xb)^\top\B^\top\}
-2\tr\{(\X-\Xb)(\X-\bTh)^\top \bSi^{-1}\B\} +np]\\
=&{1\over np}E[ \tr\{\bSi^{-1}\B(\X-\Xb)(\X-\Xb)^\top\B^\top\}
-2\tr\{(n-1) \B\} +np],
\end{align*}
because $E[(\X-\Xb)\bTh^\top]=E[\X\X^\top-\Xb\X^\top] - (n-1)\bSi$.
The quantity inside the expectation is an unbiased estimator of the risk functions.
When $n-1\geq p$, the optimal $\B$ in terms of minimizing the unbiased estimator is $\Bbh=(n-1)\bSi \{(\X-\Xb)(\X-\Xb)^\top\}^{-1}= (n-1)\W^{-1}$ for $\W$ defined in (\ref{W}), and the resulting linear shrinkage estimator is $\bThh^{LS}(\Bbh) = \X - (n-1)\W^{-1}(\X-\Xb)$, which is a variant of the Efron-Morris estimator.
It is noted that the unbiased estimator of the risk function is derived without assuming the normality.
However, $\Bbh=\W^{-1}$ is likely to be instable when $n-1$ is close to $p$.

\medskip
In this paper, we consider the ridge-type linear shrinkage estimator
\begin{equation}\label{RLS}
\bThh_{a,b}^{RLS}\coloneqq \X - \Big\{ a\Rid{1}+{b\over \tr(\W)}\I\Big\}(\X-\Xb),
\end{equation}
where $\alh$ is a nonnegative function of $\W$, and $a$ and $b$ are scalars.
In the case of $b=0$, we call $\bThh_{a,0}^{RLS}$ a {\it single shrinkage estimator }, which corresponds to Efron and Morris(1972) when $\alh=0$, and $\bThh_{a,b}^{RLS}$ is called a {\it double shrinkage estimator}, corresponding to Efron and Morris(1976) when $\alh=0$.
As an example of the ridge function $\alh$, we treat the two cases
\begin{equation}
\label{alh}
\alh = c \quad {\rm and}\quad \alh=c\, \tr \W
\end{equation}
where $c$ is a positive constant chosen suitably.

\medskip
One of features of the ridge-type linear shrinkage estimators (\ref{RLS}) is to be interpreted as the singular value shrinkage. 
In fact, in the case of $n-1\ge p$, let $\bSi^{-1/2}(\X-\Xb)=\U_s\La\V_s$ where $\U_s$ is a $p\times p$ orthogonal matrix, $\V_s$ is a $n\times n$ orthogonal matrix,  $\La=\{\diag(\si_1,\ldots,\si_p)\O_{p,n-p}\}$ and $\si_1\ge \cdots\ge \si_p$ are singular values of $\bSi^{-1/2}(\X-\Xb)$. 
Then the ridge-type linear shrinkage estimator can be expressed through the singular value decomposition as 
$$\bThh_{a,b}^{RLS}=\Xb + \bSi^{1/2}\U_s\left\{ \diag\left(\si_1\left(1-{a\over \si_1^2+\alh}-{b\over \sum \si_i^2}\right),\cdots,\si_p\left(1-{a\over \si_p^2+\alh}-{b\over \sum \si_i^2}\right)\right)\O_{p,n-p}\right\}\V_s.$$
Similarly, in the case of $p>n-1$, let $\bSi^{-1/2}(\X-\Xb)=\U_s\La\V_s$ where $\U_s$ is a $p\times p$ orthogonal matrix, $\V_s$ is a $n\times n$ orthogonal matrix,  $\La=\{\diag(\si_1,\ldots,\si_{n-1},0)\O_{n,p-n}\}^\top$ and $\si_1\ge \cdots\ge \si_{n-1}$ are singular values of $\bSi^{-1/2}(\X-\Xb)$. 
Then,
$$\bThh_{a,b}^{RLS}=\Xb + \bSi^{1/2}\U_s\left\{ \diag\left(\si_1\left(1-{a\over \si_1^2+\alh}-{b\over \sum \si_i^2}\right),\cdots,\si_{n-1}\left(1-{a\over \si_{n-1}^2+\alh}-{b\over \sum \si_i^2}\right),0\right)\O_{n,p-n}\right\}^\top\V_s.$$
Estimators treated by Efron and Morris(1972), Stein(1973) and Matsuda and Komaki(2015) are also singular value shrinkage estimators.
In the case that true mean is low rank, the singular value shrinkage method was recently studied by Candes et al(2013), Josse and Sardy(2016), Gavish and Donoho(2014, 17), among many others.

\medskip
Another feature of the ridge-type linear shrinkage estimators is that it is available for any $n$ and $p$ as long as $n\geq 2$.
Using the following lemma, which can be verified by singular value decomposition, we do not have to deal the two cases of $n-1\geq p$ and $p>n-1$ separately.

\begin{lem}\label{lem:inv}
For any $p\times (n-1)$ matrix $\A$ of which the rank is $\min\{p, n-1\}$, and any positive constant $\al$,
$$
(\A\A^\top + \al\I)^{-1}\A
=\A(\A^\top\A+\al\I)^{-1}.
$$
When $p>n-1$,
$$
(\A\A^\top)^+\A
=\A(\A^\top\A)^{-1},
$$
where $\B^+$ denotes the Moore-Penrose generalized inverse of matrix $\B$.
\end{lem}

Lemma \ref{lem:inv} allows us to extend the Efron-Morris estimators (\ref{em1}) and (\ref{em2}) to 
\begin{align}
\bThh^{em1}=&\X - (|n-p-1|-1)\W^{+}(\X-\Xb),\label{eem1}
\\
\bThh^{em2}=&\X- \Big\{{(|n-p-1|-1)}\W^{+}+{b_0 \over \tr[\W]}\I\Big\}(\X-\Xb),\label{eem2}
\end{align}
where $b_0=\min\{p^2+p-2, (n-1)^2+(n-1)-2\}$ and
$$
\W=(\X-\Xb)(\X-\Xb)^\top \bSi^{-1}.
$$
The extended Efron-Morris estimators are minimax when $|n-p-1|\geq 2$, while, as described below, the minimaxity of the ridge-type linear shrinkage estimators is guaranteed even when $p$ is close to $n-1$.

\subsection{Unbiased risk estimation and minimaxity in the case of known weights}

We here treat the case of known weights $a$ and $b$ in the ridge-type linear shrinkage estimator $\bThh_{a,b}^{RLS}$ in (\ref{RLS}) and derive conditions on $a$ and $b$ for the minimaxity under the quadratic loss.

\medskip
For simplicity, hereafter, we set $\bSi=\I_p$, because the estimator $\bThh_{a,b}^{RLS}$ is written as a function of $\bSi^{-1/2}\X$ and $\bSi^{-1/2}\X\sim \Nc(\bSi^{-1/2}\bTh, \I_p\otimes \I_n)$.
The risk difference of $\bThh_{a,b}^{RLS}$ and $\X$ is expressed as
\begin{align}
\De_{a,b}=&np\{ E[L(\bTh, \bThh_{ab}^{RLS})]- E[L(\bTh, \X)]\}
\non\\
=&
E[\tr\{ (a^2\V^2 + 2 ab\V u+b^2u^2\I)\W\}]-2E[\tr\{(a\V+bu\I)(\X-\Xb)(\X-\bTh)^\top\}],
\label{R0}
\end{align}
where $u=1/\tr(\W)$ and 
$$
\V=(\W+\alh\I)^{-1}.
$$
Using the Stein identity, we can derive an unbiased estimator of the risk difference, where the proof is given in the Appendix. 

\begin{thm}\label{thm:URE0}
Assume that $\alh=c$ or $\alh=c\tr(\W)$.
Then the unbiased estimator of the risk difference $\De$ is
\begin{align}
\Deh_{a,b}=& \tr(\V^2\W)a^2 + 2{\tr(\V\W)\over \tr(\W)}ab + {b^2\over \tr(\W)} -2(n-p-1)\tr(\V)a
\non\\
&- 2 {(n-1)pb \over \tr(\W)} - 2\alh(\tr\V)^2a + 2(2c_0+1)\tr(\V^2\W)a+ 4{b\over \tr(\W)},
\label{URE0}
\end{align}
where 
\begin{equation}\label{c0}
c_0 = \left\{\begin{array}{ll} 0 & {\rm for}\ \alh=c\\ c &{\rm for}\ \alh=c\tr(\W).\end{array}\right.
\end{equation}
\end{thm}

It is useful for expressing $\Deh_{a,b}$ in a canonical form.
Let  $\O_1$ be an $n\times (n-1)$ matrix  such that $\O_1^\top\O_1=\I_{n-1}$ and $\O_1\O_1^\top=\I-n^{-1}\text{\bf 1}_n\text{\bf 1}_n^\top$. 
Let $\Y=\X\O_1$ and $\bxi=\bTh\O_1$. 
Then,  it follows that $\X-\Xb=\X\O_1\O_1^\top$ and $\Y\sim \Nc_{p, n-1}(\bxi, \I_p\otimes \I_{n-1})$.
In the case of $n-1\geq p$, let $\S=\Y\Y^\top$ and $\U=(\S+\al\I)^{-1}$.
Since $u=1/\tr(\W)=1/\tr(\S)$, we have
\begin{align*}
\Deh_{a,b}
=& \tr(\U^2\S)a^2 + 2\tr(\U\S)uab + ub^2 -2(n-p-1)\tr(\U)a
\non\\
&- 2 (n-1)pub - 2\alh(\tr\U)^2a + 2(2c_0+1)\tr(\U^2\S)a+ 4ub.
\end{align*}
In the case of $p>n-1$, let $\S'=\Y^\top\Y$ and $\U'=(\S'+\al\I)^{-1}$.
Then, $u=1/\tr(\W)=1/\tr(\S')$ and $\Deh_{a,b}$ is expressed as
\begin{align}
\Deh_{a,b}
=& \tr(\U'^2\S')a^2 + 2\tr(\U'\S')uab + ub^2 -2(p-n+1)\tr(\U')a
\non\\
&- 2 (n-1)pub - 2\alh(\tr\U')^2a + 2(2c_0+1)\tr(\U'^2\S')a+ 4ub,
\label{URE0-2}
\end{align}

We first treat the single shrinkage estimators which corresponds to $\bThh_{a,0}^{RLS}$ for $b=0$.
The conditions on $a$ for minimaxity of $\bThh_{a,0}^{RLS}$ are provided from Theorem \ref{thm:URE0}.

\begin{prp}\label{prp:1}
For $\alh$ and $c_0$ given in $(\ref{alh})$ and $(\ref{c0})$, $\bThh_{a,0}^{RLS}$ is minimax if $a$ satisfies the condition
$$
0<a \leq 2[ |n-p-1|-1 + \{(n-1)p-2\}c_0].
$$
\end{prp}

{\it Proof}.\ \ 
Consider the case of $n-1\geq p$.
From (\ref{URE0}), the unbiased estimator of the risk difference is 
\begin{equation}
\Deh_{a,0}= \tr(\V^2\W)a \Big\{ a   -2 {(n-p-1)\tr(\V)+\alh(\tr\V)^2 \over  \tr(\V^2\W)} + 2(2c_0+ 1)\Big\}
\label{Dea0}
\end{equation}
In the case of $\alh=c$, since $\tr(\V^2\W)\leq \tr(\V)$, $\Deh_{a,0}$ is evaluated as
$$
\Deh_{a,0}\leq \tr(\V^2\W)a \Big\{ a   -2 \{(n-p-1)-1\} -2 {c (\tr\V)^2\over \tr(\V^2\W)}\Big\},
$$
which gives the sufficient condition $0<a\leq 2 \{(n-p-1)-1\}$.
In the case of $\alh=c\tr\W$, from Lemma \ref{lem:ineq} in appendix, note that $\tr(\V^2\W)\leq \tr(\V)/(1+pc)$ and $\tr(\W)\tr(\V)\geq p^2/(1+pc)$.
Then from (\ref{Dea0}), we have
$$
\Deh_{a,0}\leq \tr(\V^2\W)a \Big\{ a   -2 {(n-p-1)\tr(\V)+c\tr(\V) p^2/(1+pc) \over  \tr(\V)/(1+pc)} + 2(2c+ 1)\Big\},
$$
which yields the sufficient condition $0<a\leq 2 [ \{(n-p-1)-1\} + \{(n-1)p-2\}c]$.
In the case of $p>n-1$, from the equation (\ref{URE0-2}), we get the corresponding conditions.
\hfill$\Box$

\bigskip
Proposition \ref{prp:1} implies that the ridge-type linear shrinkage estimator with $\alh=c\tr\W$ is minimax for $c>1/\{(n-1)p-2\}$ and $(n-1)p>2$ even when $n-1=p$.

\medskip
We next consider the double shrinkage estimators which corresponds to $\bThh_{a,b}^{RLS}$ for $b\not=0$.
The condition on $a$ and $b$ for minimaxity of $\bThh_{a,b}^{RLS}$ is provided from Theorem \ref{thm:URE0}.

\begin{prp}\label{prp:2}
For $\alh$ and $c_0$ given in $(\ref{alh})$ and $(\ref{c0})$, the double shrinkage estimator $\bThh_{a,b}^{RLS}$ dominates the single shrinkage estimator $\bThh_{a,0}^{RLS}$ if 
$$
0<b\leq2(n-1)p - 4 - 2a\min(p, n-1)/\{1+c_0\min(p, n-1)\}.
$$
Combining this result and Proposition \ref{prp:1} gives the condition for minimaxity, namely $\bThh_{a,b}^{RLS}$ is minimax provided 
$$
0<a \leq 2[ |n-p-1|-1 + \{(n-1)p-2\}c_0]\ {\rm and}\ 0<b\leq 2(n-1)p - 4 - 2{a\min(p, n-1)\over 1+c_0\min(p, n-1)}.
$$
When $a=|n-p-1|-1 + \{(n-1)p-2\}c_0$, the condition on $b$ for the minimaxity is 
$$
0<b\leq 2{ (n-1)p-2-min(p, n-1)(|n-p-1|-1) \over 1+c_0\min(p, n-1)}.
$$
\end{prp}

{\it Proof}.\ \ 
From (\ref{URE0}) and (\ref{Dea0}), the unbiased estimator of the risk difference is 
$$
\Deh_{a,b}= \Deh_{a,0} +{b\{b-2(n-1)p+4+2\tr(\V\W)a\}\over \tr(\W)}.
$$
In the case of $n-1\geq p$, note that $\tr(\V\W)\leq p$.
Then, $\Deh_{a,b}\leq \Deh_{a,0}$ if $b\leq2(n-1)p - 4 - 2pa$.
When $\alh=c\tr(\W)$, from Lemma \ref{lem:ineq}, we have $\tr(\V\W)\leq p/(1+cp)$, so that $\Deh_{a,b}\leq \Deh_{a,0}$ if $b\leq2(n-1)p - 4 - 2pa/(1+cp)$.
Substituting $a=|n-p-1|-1 + \{(n-1)p-2\}c_0$, one gets the sufficient conditions.
\hfill$\Box$

\section{Proposal of Shrinkage Weights and Minimaxity in the Case of Estimated Weights}
\label{sec:PWMstein}

The case of known $a$ and $b$ is treated in the previous section.
However, $a$ and $b$ are shrinkage weights to control the extent of shrinkage, and we want to determine them from data.
To this end, we obtain the shrinkage weights in terms of minimizing the unbiased risk estimator  (\ref{URE0}).

\medskip
We first consider the single shrinkage estimator $\bThh_{a,0}^{RLS}$ for $b=0$.
By minimizing (\ref{URE0}) with respect to $a$ when $b=0$, we get 
\begin{equation}\label{ash}
\ah_S = {(n-p-1)\tr(\V) + \alh (\tr\V)^2 \over \tr(\V^2\W)} - (2c_0+1).
\end{equation}
Corresponding to the two cases of $\alh$, we get the ridge-type single linear shrinkage estimators $\bThh^{S1}=\bThh_{\ah_S,0}^{RLS}$ for $\alh=c$ and $\bThh^{S2}=\bThh_{\ah_S,0}^{RLS}$ for $\alh=c\tr(\W)$. 

\medskip
We investigate the limiting values of the weight $\ah_S$ as $c\to 0$ or $\infty$. 
For $\alh=c\ {\text or}\ c\tr(\W)$ and any natural number $c_*$, it can be seen that $\alh^{c_*}\tr(\V^{c_*})\to p\as$ and $\alh^{c_*}\tr(\V^{c_*}\W)\to\tr(\W)\as\asl c\to \infty$. 
Using these limits, we can provide the following result.

\begin{prp}
\label{prp:alim}
Let $\si_i$ be singular values of $\X-\Xb$.
In the case of $|n-p-1|>1$, it holds that
\begin{align*}
\ah_S \to|n-p-1|-1 \as\asl c\to 0,
\end{align*}
and
\begin{align*}
{\ah_S\over \si^2_i+\alh}&\to \left\{\begin{array}{ll} {(n-1)p/ \tr(\W)} & {\rm for}\ \alh=c\\ \{(n-1)p-2\}/\tr(\W)\} &{\rm for}\ \alh=c\tr(\W)\end{array} \as \asl c \to \infty. \right.\\
\end{align*}
\end{prp}

Proposition \ref{prp:alim} shows that the weight approaches that of the extended Efron-Morris estimator when $c\to 0$, and as $c \to\infty$, the weight approaches that of the James-Stein estimator 
\begin{align}
\bThh^{js}=(1-b_{js}/\tr[\W])(\X-\Xb)+\Xb \label{js}
\end{align}
for ($b_{js}=(n-1)p$ or) $b_{js}=\{(n-1)p-2\}$. 

\medskip
For the ridge-type linear shrinkage estimator $\bThh_{\ah_S,0}^{RLS}$ with the estimated weight $\ah_S$, it is quite interesting to investigate the minimaxity, but harder than Proposition \ref{prp:1}.

\begin{thm}\label{thm:min}
The single shrinkage estimator $\bThh_{\ah_S,0}^{RLS}$ with $\ah_S$ substituted is minimax if either of the following conditions is satisfied in the case of $n-1\geq p$:\\

{\rm (i)}\ $\alh=c$ and $n-p-1\geq 10$ \\

{\rm (ii)}\  $\alh=c\tr\W$ and
\begin{equation}\label{minc1}
n-p-1 \geq { (-p^2+12)c^2 + (-p^2+8p+14+4/p)c+14\over (1+c)(1+cp)},
\end{equation}
or equivalently, for $A_0=n-p-1$,
\begin{equation}\label{minc2}
(pA_0+p^2-12)c^2 + \{(p+1)A_0+p^2-8p-14-4/p\}c + A_0-14\geq 0.
\end{equation}

{\rm (iii)}\  $\alh=(1/\min(n-1,p))\tr\W$ and 
\begin{equation}\label{minc3}
n-p-1 \geq {-p^3+21p^2+14p+16\over 2p(p+1)}.
\end{equation}

In the case of $p>n-1$, the conditions for the mimimaxity are given by replacing $n-1$ and $p$ in the above conditions with $p$ and $n-1$, respectively.
\end{thm}

The proof is given in the Appendix.
The conditions (ii) and (iii) include the case of $n-1=p$, while even the extended Efron-Morris estimator does not achieve the minimaxity. For example, when $p=4$, the numerator of the right hand side of (\ref{minc1}) is $-4c^2+31c+14$, and it is negative for $c\geq 8.178$. In that case, the minimaxity is achieved even when $n-1=p$. We can see that when $p$ increases, the lower bound of $c$ for which the proposed estimator holds minimaxity approaches 0.

The right hand side of condition (iii) is negative for $p\geq 22$. When $p\geq 22$, minimaxity is always achieved. The value of right hand side of condition (iii) is 12.5 when $p=1$ and it is 10 when $p=2$. The value decreases toward 0 as p increases toward 22.

\medskip
We can consider a positive-part estimator like Efron and Morris (1972) and Tsukuma (2010). 
Let $\bThh^{RLS+}_{\ah_S,0}$ be the estimator which is obtained by replacing $\si_i(1-\ah_S/(\si_i^2+\alh))$ with $\max(0,\si_i(1-\ah_S/(\si_i^2+\alh)))$ in $\bThh^{RLS}_{\ah_S,0}$.
Then the positive-part estimator dominates original estimator. The proof is the same way as Tsukuma (2010).

\begin{cor}\label{cor:min}
The positive-part single shrinkage estimator $\bThh_{\ah_S,0}^{RLS+}$ dominates the single shrinkage estimator $\bThh_{\ah_S,0}^{RLS}$. The estimator $\bThh_{\ah_S,0}^{RLS+}$ is also minimax under conditions in Theorem \ref{thm:min}.
\end{cor}

The positive-part ridge-type linear shrinkage estimator for $\alh=c\tr(\W)$ is denoted by $\bThh^{S2+}=\bThh_{\ah_S,0}^{RLS+}$, and it will be used for simulation study in Section \ref{sec:sim} as the proposal single shrinkage estimator.

\medskip
We next derive the weight estimator $\bh_S$ in terms of minimizing $\Deh_{\ah_S,b}$ with respect to $b$ for given $\ah_S$. 
It is noted that 
$$
\Deh_{\ah_S,b}=-\tr[\V^2\W]\ah_S^2-4\sum{\ell_i\over \ell_i+\alh}{\partial \ah_S\over \partial \ell_i}+2\tr[\V\W]u\ah_Sb+ub^2-2(n-1)pub+4ub,
$$
which is minimized at
\begin{equation}\label{bh}
\bh_S = \{(n-1)p-2\}-\tr(\V\W)\ah_S.
\end{equation}
The ridge-type double shrinkage estimators are denoted by $\bThh^{D1}=\bThh_{\ah_S,\bh_S}^{RLS}$ for $\alh=c$ and $\bThh^{D2}=\bThh_{\ah_S,\bh_S}^{RLS}$ for $\alh=c\tr(\W)$. 

\begin{prp}\label{prp:blim}
Corresponding to Proposition \ref{prp:alim}, $\bh_S$ converges to the following limits:\\

$\bh_S\to \min\{(n-1)^2+(n-1)-2, p^2+p-2\} \as \asl c\to 0$.\\

$\bh_S\to\left\{\begin{array}{cl} -2 & {\rm for}\ \alh=c\\ 0 &{\rm for}\ \alh=c\tr[\W]\end{array}\right. \as \asl c\to\infty$.
\end{prp}

From Propositions \ref{prp:alim} and \ref{prp:blim},  it is seen that the proposed double shrinkage estimator approaches the extended Efron-Morris estimator $\bThh^{em2}$ in $(\ref{eem2})$ when $c\to 0$, and as $c\to \infty$, the proposed estimator approaches the James-Stein estimator $\bThh^{js}=(1-b_{js}/\tr[\W])(\X-\Xb)+\Xb$ for $b_{js}=\{(n-1)p-2\}$.

We can show that double shrinkage estimator also has minimaxity.

\begin{thm}\label{thm:dmin}
The double shrinkage estimator $\bThh_{\ah_S,\bh_S}^{RLS}$ is minimax if either of the following conditions is satisfied in the case of $n-1\geq p$:\\

{\rm (i)}\  $\alh=c\tr\W$ and
\begin{equation}\label{minc1d}
n-p-1 \geq { (-p^2+16-4/p)c^2 + (-p^2+8p+18)c+14\over (1+c)(1+cp)},
\end{equation}
or equivalently, for $A_0=n-p-1$,
\begin{equation}\label{minc2d}
(pA_0+p^2-16+4/p)c^2 + \{(p+1)A_0+p^2-8p-18\}c + A_0-14\geq 0.
\end{equation}

{\rm (ii)}\  $\alh=(1/\min(n-1,p))\tr\W$ and 
\begin{equation}\label{minc3d}
n-p-1 \geq {-p^4+21p^3+18p^2+16p-4\over 2p^2(p+1)}.
\end{equation}

In the case of $p>n-1$, the conditions for the mimimaxity are given by replacing $n-1$ and $p$ in the above conditions with $p$ and $n-1$, respectively.
\end{thm}

The proof is given in the Appendix.
The sufficient conditions include the case of $n-1=p$, while even the extended Efron-Morris estimator does not achieve the minimaxity. For example, when $p=4$, the numerator of the right hand side of (\ref{minc1d}) is $-c^2+34c+14$, and it is negative for $c\geq 34.407$. In that case, the minimaxity is achieved even when $n-1=p$. We can see that when $p$ increases, the lower bound of $c$ for which the proposed estimator holds minimaxity approaches 0.

The right hand side of condition (ii) is negative for $p\geq 22$. When $p\geq 22$, minimaxity is always achieved. The value of right hand side of condition (ii) is 12.5 when $p=1$ and it is 10.5 when $p=2$. The value decreases toward 0 as p increases toward 22.

The positive-part double shrinkage estimator $\bThh^{RLS+}_{\ah_S,\bh_S}$ is given by replacing $\si_i(1-\ah_S/(\si_i^2+\alh)-\bh_S/\sum\si_i^2)$ with $\max(0,\si_i(1-\ah_S/(\si_i^2+\alh)-\bh_S/\sum\si_i^2))$ in $\bThh^{RLS}_{\ah_S,\bh_S}$. 
Then the positive-part estimator dominates the original one.

\begin{cor}\label{cor:dmin}
The positive-part double shrinkage estimator $\bThh_{\ah_S,\bh_S}^{RLS+}$ dominates the double shrinkage estimator $\bThh_{\ah_S,\bh_S}^{RLS}$. 
The estimator $\bThh_{\ah_S,\bh_S}^{RLS+}$ is also minimax under conditions in Theorem \ref{thm:dmin}.
\end{cor}

The positive-part ridge-type linear shrinkage estimator for $\alh=c\tr(\W)$ is denoted by $\bThh^{D2+}=\bThh_{\ah_S,\bh_S}^{RLS+}$ and it is used for simulation study in Section \ref{sec:sim} as the proposal double shrinkage estimator.

\section{A Relation with the Optimal Weights Derived from the Random Matrix Theory}
\label{sec:RMT}

We here derive the optimal weights based on Bayesian framework. Then we show that our proposed weights in previous section has consistency using the random matrix theory.
Assume that $\bSi=\I$ without any loss of generality.
Consider again the prior distribution $\bth_i \mid \bmu \sim \Nc_p(\bmu, \bLa)$ and $\bmu\sim {\rm Uniform\  on\ } \Re^p$.
Let $\Y=\X\O_1$ where $\O_1$ be an $n\times (n-1)$ matrix such that $\O_1^\top\O_1=\I_{n-1}$ and $\O_1\O_1^\top=\I-n^{-1}\text{\bf 1}_n\text{\bf 1}_n^\top$.
Then, the Bayes estimator is $\bThh^{Bays}=\X-\bPsi^{-1}(\X-\Xb)$ for $\bPsi=\I+\bLa$, and the marginal distribution of $\Y$ is $\Nc(\zero,  \bPsi\otimes \I_{n-1})$.
When $\bPsi$ is estimated by $\bPsih$, the risk function of $\bThh^{LS}(\bPsih^{-1})=\X-\bPsih^{-1}(\X-\Xb)$ relative to the quadratic loss (\ref{loss}) can be expressed as
\begin{align*}
{1\over np}E[L(\bTh,\bThh^{LS}(\bPsih^{-1}))]=&{1\over np}E[E[\tr(\bTh-\bThh^{Bayes})^\top(\bTh-\bThh^{Bayes})\mid \X] ]\\
&
+ {1\over np}E[\tr(\bPsih^{-1}-\bPsi^{-1})^2(\X-\Xb)(\X-\Xb)^\top],
\end{align*}
where the expectation is taken with respect to the joint density of $(\X, \bTh)$.
Let $\W=(\X-\Xb)(\X-\Xb)^\top=\Y\Y^\top$.
Thus, the problem is reduced to the estimation of the precision matrix $\bPsi^{-1}$ in the model $\Y\sim\Nc(\zero,  \bPsi\otimes \I_{n-1})$ relative to the loss function $L_{EM}(\bPsi, \bde)=(np)^{-1}\tr\{(\bde-\bPsi^{-1})^2\W\}$, which was treated by Efron and Morris (1976).

\medskip
Using the random matrix theory, we can derive the optimal $a$ in the single shrinkage estimator
$$
\bde_a = a(\W+\alh\I_p)^{-1}
$$
under the loss $L_{EM}(\bPsi,\bde)$.
When $n-1\geq p$, the matrix $\W$ is non-singular, and the loss function of $\bde_a$ is
\begin{equation*}
L_{EM}(\bPsi,\bde_a)={1\over np}\Big[a^2 \tr\{ (\W+\alh\I_p)^{-2}\W\}-2a\tr\{(\W+\alh\I_p)^{-1}\bPsi^{-1}\W\} + \tr\{\bPsi^{-2}\W\}\Big].
\end{equation*}
The optimal $a$ is
\begin{equation}
a^*={\tr[\bPsi^{-1}]-\alh\tr[\Rid{1}\bPsi^{-1}]\over \trRS{2}}. \label{ahoptR1}
\end{equation}
It is interesting to investigate the relation between the optimal $a^*$ and the estimator $\ah_S$ given in (\ref{ash}).
We here show that the estimator $\ah_S$ converges to the optimal $a^*$ by using random matrix theory.
For this purpose, we assume the following conditions.

\smallskip
(A1) $n/p \to c \asl n\to\infty$, where $c$ is constant and $0<c<\infty$ and $c\neq 1$.

\smallskip
(A2) The eigenvalues of $\bPsi$ is uniformly bounded that is there are positive constants $c_1, c_2$ satisfying $c_1\leq \min\{\tau_i\}\leq\max\{\tau_i\}\leq c_2$ where $\tau_1, \ldots, \tau_p$ are eigenvalues of $\bPsi$. It is also assumed that $H_p(\tau)=p^{-1}\sum_{j=1}^p 1_{[\tau_j,\infty)}(\tau)$ converges to nonrandom probability distribution $H(\tau)$.

\smallskip
(A3) $\alh/n\to \al\asl n\to\infty\as$, where $\al$ is positive.

\begin{thm}
\label{thm:1}
Under conditions {\rm(A1)-(A3)}, it holds that $|\ah_S-a^*|/n\to 0\asl n\to\infty\as$
\end{thm}

We next consider the double shrinkage estimator
$$
\bde_{\ah_S,b} = \ah_S(\W+\alh\I_p)^{-1} + {b\over \tr(\W)}\I.
$$
which has the loss
\begin{align*}
L_{EM}(\bPsi,\bde_{\ah_S,b})=& {1\over np}\Big[{\ah_S}^2\tr(\V^2\W)+2\ah_Sb\tr(\V\W)/\tr(\W)+b^2/\tr(\W)\\
&-2\ah_S\tr(\bPsi^{-1})+2\ah_S\alh\tr(\V\bPsi^{-1})-2b\tr(\W\bPsi^{-1})/\tr(\W)+\tr(\W\bPsi^{-2})\Big],
\end{align*}
for $\V=(\W+\alh\I_p)^{-1}$.
By differentiating the loss with respect to $b$, the optimal $b$ is given by
\begin{align}
b^*=-\ah_S\tr(\V\W)+\tr(\W\bPsi^{-1}). \label{bhoptR2}
\end{align}

\begin{thm}
\label{thm:2}
Under conditions {\rm(A1)-(A3)}, it holds that $|\bh_S/\tr\W-b^*/\tr\W|\to 0\asl n\to\infty\as$
\end{thm}

\section{Simulation Study}
\label{sec:sim}
We now investigate the numerical performance of the suggested ridge-type linear shrinkage estimators and compare the performance with  the Efron-Morris estimator, the James-Stein estimator and the singular value shrinkage estimator proposed by Gavish and Donoho (2017).

\medskip
The six estimators we compare are $\bThh^{S2+}$, $\bThh^{D2+}$, $\bThh^{em+}$, $\bThh^{em2+}$, $\bThh^{js+}$ and $\bThh^{gd}$, denoted by S2+, D2+, em+, em2+, js+ and gd.
These estimators are defined as follows:

$\bThh^{S2+}$ is $\bThh^{RLS+}_{\ah_{S},0}$ with $\alh=\tr\W/\min(n-1,p)$.

$\bThh^{D2+}$ is $\bThh^{RLS}_{\ah_S,\bh_S}$ with $\alh=\tr\W/\min(n-1,p)$. 

$\bThh^{em+}$ and $\bThh^{em2+}$ are positive-part versions of the extended Efron-Morris estimator given in (\ref{eem1}) and (\ref{eem2}).

$\bThh^{js+}$ is a positive-part version of the James-Stein estimator given in (\ref{js}).

$\bThh^{gd}$ is the singular value shrinkage estimator suggested by Gavish and Donoho (2017), namely,
\begin{align*}
\bThh^{gd}=&\Xb+\U_s \diag\Big\{\sqrt{((\si_1^2/(n-1)-p/(n-1)-1)^2-4p/(n-1))_+}/(\si_1/\sqrt{n-1}),\ldots,\\
& \hspace{3cm} \sqrt{((\si_p^2/(n-1)-p/(n-1)-1)^2-4p/(n-1))_+}/(\si_p/\sqrt{n-1})\Big\}\V_s,
\end{align*}
where $a_+=\max(a, 0)$.
This estimator is available for $n-1\geq p$, while for $p>n-1$, we can calculate the corresponding estimator for $\X^\top$.
The estimator works well in the situation that the true mean matrix is low rank.

\medskip
The simulation experiments are conducted by four cases of data generated from a normal distribution and two cases from non-normal distributions.
We first consider the case that $\X$ is generated from the normal distribution with mean $\bTh$ and the identity covariance matrix.
When $s_1,\ldots,s_{\min(n,p)}$ are singular values of $\bTh$, we consider the two setups for $s_i$.
The first setup is that $s_i=10+10\times (i-1)/\{(\min(n,p)/5)-1\}$ for $i=1,\ldots, \min(n,p)/5$ and the others are $s_i=10^{q}$, where the power $q$ controls the situation with/without low rank. 
In this experiment, we take the values $q=-1, 1/2$, and the corresponding risk averages of simulated values based on 5,000 replications are given in Table \ref{sim1-3} and Table \ref{sim1-2}, respectively.
The second setup is the cases of larger singular values, namely, $s_i=100+100\times (i-1)/\{(\min(n,p)/10)-1\}$ for $i=1,\ldots, \min(n,p)/10$ and $s_i=\min(n,p)^{q}$. 
Risk averages of simulated values based on 5,000 replications for $q=-1, 1/2$ are given in Table \ref{sim2-3} and  Table \ref{sim2-2}), respectively. 

\medskip
We next investigate the performance for non-normal distributions with mean $\bTh$ and the identity covariance matrix, where singular values are given by $s_i=10+10\times (i-1)/\{(\min(n,p)/5)-1\}$ for $i=1,\ldots, \min(n,p)/5$ and the others are $s_i=10^{-1}$.
The non-normal distributions we examine are the t-distribution with 3 degrees of freedom and the  chi-square distribution with 2 degrees of freedom, and the corresponding risk averages of simulated values based on 5,000 replications are given in Table \ref{simt1-3} and Table \ref{simchi1-3}, respectively.

\medskip
The six tables report the averages of simulated risk values for the above six estimators, where the setups of $(n,p)$ are $(n,p)=(100,10)$, $(100,30)$, $(100,80)$, $(100,101)$, $(10,100)$, $(30,100)$ and $(80,100)$.
It is noted that the risk of the estimator $\X$ is one.
From these numerical investigations, we observe that the extended Efron-Morris estimators based on the Moore-Penrose generalized inverse do not work in the case of $n-1=p$, while the suggested ridge-type linear shrinkage estimators are stable and the best.
In the situation with low rank, the Gavish-Donoho estimator is the best and the suggested estimators perform reasonably well.
In the situation without low rank, the suggested estimators work better than the Efron-Morris and James-Stein estimators.

\medskip
The performances in the low rank situation are reported in Tables \ref{sim1-3} and \ref{sim2-3} for the normal distribution and in Tables \ref{simt1-3} and \ref{simchi1-3} for the t- and chi-square distributions.
Since the true mean has low rank, the situation is advantageous for the Gavish-Donoho estimator $\bThh^{gd}$. 
In fact, the performance of the Gavish-Donoho estimator is the best for the normal distribution.
However, the suggested estimators are better than $\bThh^{gd}$ for the t-distribution and comparable to $\bThh^{gd}$ for the chi-square distribution.
Especially, as seen from Table \ref{sim2-3}, the Efron-Morris and James-Stein estimators do not work well, but the suggested double shrinkage estimator $\bThh^{D2+}$ performs better.

\medskip
Tables \ref{sim1-2} and \ref{sim2-2} examine the situations that the true mean does not have low rank.
From Table \ref{sim1-2}, it is revealed that the suggested estimators $\bThh^{S2+}$ and $\bThh^{D2+}$ are comparable to $\bThh^{gd}$ and better than the Efron-Morris estimators.
When $n$ is close to $p$,  the Efron-Morris estimators does not work well, and the suggested estimators are good and better than $\bThh^{gd}$.
In the situation reported in Table \ref{sim2-2}, the risk gain of the James-Stein estimator is small.
When $n$ is close to $p$, the suggested estimators perform better than the other ones. We verified by another setting Gavish-Donoho estimator does not have minimaxity, that is, the simulation risk of $\bThh^{gd}$ could be grater than 1.

\begin{table}[H]
\caption{\small Simulated Risk Values of the Six Estimators where $s_i=10+10\times (i-1)/\{(\min(n,p)/5)-1\}$ for $i=1,\ldots, \min(n,p)/5$ and the others are $s_i=10^{-1}$.
}
\centering
$
{\renewcommand\arraystretch{1.1}
\begin{array}{c@{\hspace{1mm}}
              r@{\hspace{1mm}}
              r@{\hspace{1mm}}
              r@{\hspace{1mm}}
              r@{\hspace{1mm}}
              r@{\hspace{1mm}}
              r@{\hspace{1mm}}
              r
             }
\text{$(n, p)$} &\text{S2+}&\text{D2+}&\text{em+}&\text{em2+}&\text{js+}&\text{gd}\\
\hline
\text{$(100, 10)$}
&0.181
&0.181
&0.204
&0.184
&0.368
&0.165
\\
\text{$(100, 30)$} 
&0.217
&0.217
&0.342
&0.237
&0.333
&0.216
\\
\text{$(100, 80)$} 
&0.273
&0.263
&0.775
&0.294
&0.324
&0.299
\\
\text{$(101, 100)$} 
&0.290
&0.276
&34.5
&33.8
&0.325
&0.321
\\
\text{$(10, 100)$} 
&0.258
&0.258
&0.274
&0.259
&0.375
&0.166
\\ 
\text{$(30, 100)$} 
&0.237
&0.237
&0.351
&0.255
&0.346
&0.216
\\
\text{$(80, 100)$}
&0.275
&0.266
&0.758
&0.296
&0.328
&0.302
\\
\hline
\end{array}
}
$
\label{sim1-3}
\end{table}
\normalsize

\begin{table}[H]
\caption{\small Simulated Risk Values of the Six Estimators  where $s_i=10+10\times (i-1)/\{(\min(n,p)/5)-1\}$ for $i=1,\ldots, \min(n,p)/5$ and the others are $s_i=10^{1/2}$.
}
\centering
$
{\renewcommand\arraystretch{1.1}
\begin{array}{c@{\hspace{1mm}}
              r@{\hspace{1mm}}
              r@{\hspace{1mm}}
              r@{\hspace{1mm}}
              r@{\hspace{1mm}}
              r@{\hspace{1mm}}
              r@{\hspace{1mm}}
              r
             }
\text{$(n, p)$} &\text{S2+}&\text{D2+}&\text{em+}&\text{em2+}&\text{js+}&\text{gd}\\
\hline
\text{$(100, 10)$}
&0.250
&0.250
&0.274
&0.256
&0.400
&0.245
\\
\text{$(100, 30)$} 
&0.287
&0.286
&0.403
&0.304
&0.368
&0.295
\\
\text{$(100, 80)$} 
&0.343
&0.331
&0.799
&0.342
&0.359
&0.377
\\
\text{$(101, 100)$} 
&0.359
&0.343
&143
&143
&0.360
&0.397
\\
\text{$(10, 100)$} 
&0.319
&0.319
&0.337
&0.323
&0.409
&0.245
\\ 
\text{$(30, 100)$} 
&0.304
&0.303
&0.411
&0.320
&0.380
&0.294
\\
\text{$(80, 100)$}
&0.345
&0.333
&0.783
&0.344
&0.363
&0.379
\\
\hline
\end{array}
}
$
\label{sim1-2}
\end{table}

\begin{table}[H]
\caption{\small Simulated Risk Values of the Six Estimators  where $s_i=100+100\times (i-1)/\{(\min(n,p)/10)-1\}$ for $i=1,\ldots, \min(n,p)/10$ and the others are $s_i=10^{-1}$.
}
\centering
$
{\renewcommand\arraystretch{1.1}
\begin{array}{c@{\hspace{1mm}}
              r@{\hspace{1mm}}
              r@{\hspace{1mm}}
              r@{\hspace{1mm}}
              r@{\hspace{1mm}}
              r@{\hspace{1mm}}
              r@{\hspace{1mm}}
              r
             }
\text{$(n, p)$} &\text{S2+}&\text{D2+}&\text{em+}&\text{em2+}&\text{js+}&\text{gd}\\
\hline
\text{$(100, 10)$}
&0.203
&0.149
&0.172
&0.169
&0.914
&0.118
\\
\text{$(100, 30)$} 
&0.326
&0.241
&0.310
&0.305
&0.961
&0.136
\\
\text{$(100, 80)$} 
&0.352
&0.269
&0.765
&0.736
&0.959
&0.181
\\
\text{$(101, 100)$} 
&0.367
&0.285
&99.0
&99.0
&0.960
&0.198
\\
\text{$(10, 100)$} 
&0.297
&0.243
&0.250
&0.248
&0.918
&0.120
\\ 
\text{$(30, 100)$} 
&0.346
&0.261
&0.320
&0.315
&0.961
&0.137
\\
\text{$(80, 100)$}
&0.359
&0.275
&0.746
&0.719
&0.960
&0.183
\\
\hline
\end{array}
}
$
\label{sim2-3}
\end{table}
\normalsize

\begin{table}[H]
\caption{\small Simulated Risk Values of the Six Estimators  where $s_i=100+100\times (i-1)/\{(\min(n,p)/10)-1\}$ for $i=1,\ldots, \min(n,p)/10$ and the others are $s_i=10^{1/2}$.
}
\centering
$
{\renewcommand\arraystretch{1.1}
\begin{array}{c@{\hspace{1mm}}
              r@{\hspace{1mm}}
              r@{\hspace{1mm}}
              r@{\hspace{1mm}}
              r@{\hspace{1mm}}
              r@{\hspace{1mm}}
              r@{\hspace{1mm}}
              r
             }
\text{$(n, p)$} &\text{S2+}&\text{D2+}&\text{em+}&\text{em2+}&\text{js+}&\text{gd}\\
\hline
\text{$(100, 10)$}
&0.270
&0.225
&0.250
&0.248
&0.915
&0.207
\\
\text{$(100, 30)$} 
&0.459
&0.405
&0.491
&0.485
&0.961
&0.400
\\
\text{$(100, 80)$} 
&0.624
&0.598
&0.916
&0.889
&0.960
&0.844
\\
\text{$(101, 100)$} 
&0.674
&0.654
&3.02
&2.98
&0.961
&1.00
\\
\text{$(10, 100)$} 
&0.356
&0.310
&0.321
&0.319
&0.919
&0.208
\\ 
\text{$(30, 100)$} 
&0.474
&0.420
&0.497
&0.492
&0.961
&0.398
\\
\text{$(80, 100)$}
&0.627
&0.600
&0.909
&0.883
&0.961
&0.844
\\
\hline
\end{array}
}
$
\label{sim2-2}
\end{table}

\begin{table}[H]
\caption{\small Simulated Risk Values of the Six Estimators  where $\x_i=\th_i+ t_3\sqrt{(3-2)/3}$, $s_i=10+10\times (i-1)/\{(\min(n,p)/5)-1\}$ for $i=1,\ldots, \min(n,p)/5$ and the others are $s_i=10^{-1}$.
}
\centering
$
{\renewcommand\arraystretch{1.1}
\begin{array}{c@{\hspace{1mm}}
              r@{\hspace{1mm}}
              r@{\hspace{1mm}}
              r@{\hspace{1mm}}
              r@{\hspace{1mm}}
              r@{\hspace{1mm}}
              r@{\hspace{1mm}}
              r
             }
\text{$(n, p)$} &\text{S2+}&\text{D2+}&\text{em+}&\text{em2+}&\text{js+}&\text{gd}\\
\hline
\text{$(100, 10)$}
&0.244
&0.248
&0.284
&0.261
&0.396
&0.272
\\
\text{$(100, 30)$} 
&0.274
&0.278
&0.396
&0.287
&0.354
&0.299
\\
\text{$(100, 80)$} 
&0.324
&0.321
&0.789
&0.318
&0.344
&0.376
\\
\text{$(101, 100)$} 
&0.333
&0.326
&301
&300
&0.337
&0.388
\\
\text{$(10, 100)$} 
&0.321
&0.324
&0.350
&0.333
&0.442
&0.280
\\ 
\text{$(30, 100)$} 
&0.288
&0.292
&0.399
&0.300
&0.365
&0.297
\\
\text{$(80, 100)$}
&0.321
&0.318
&0.766
&0.314
&0.341
&0.374
\\
\hline
\end{array}
}
$
\label{simt1-3}
\end{table}

\begin{table}[H]
\caption{\small Simulated Risk Values of the Six Estimators  where $\x_i=\th_i+(\chi_2^2-2)/\sqrt{2\cdot 2}$, $s_i=10+10\times (i-1)/\{(\min(n,p)/5)-1\}$ for $i=1,\ldots, \min(n,p)/5$ and the others are $s_i=10^{-1}$.
}
\centering
$
{\renewcommand\arraystretch{1.1}
\begin{array}{c@{\hspace{1mm}}
              r@{\hspace{1mm}}
              r@{\hspace{1mm}}
              r@{\hspace{1mm}}
              r@{\hspace{1mm}}
              r@{\hspace{1mm}}
              r@{\hspace{1mm}}
              r
             }
\text{$(n, p)$} &\text{S2+}&\text{D2+}&\text{em+}&\text{em2+}&\text{js+}&\text{gd}\\
\hline
\text{$(100, 10)$}
&0.190
&0.190
&0.222
&0.200
&0.370
&0.179
\\
\text{$(100, 30)$} 
&0.224
&0.224
&0.355
&0.245
&0.331
&0.224
\\
\text{$(100, 80)$} 
&0.279
&0.270
&0.778
&0.298
&0.327
&0.306
\\
\text{$(101, 100)$} 
&0.294
&0.282
&296
&295
&0.325
&0.325
\\
\text{$(10, 100)$} 
&0.277
&0.276
&0.297
&0.281
&0.417
&0.183
\\ 
\text{$(30, 100)$} 
&0.245
&0.245
&0.364
&0.265
&0.351
&0.224
\\
\text{$(80, 100)$}
&0.281
&0.272
&0.761
&0.299
&0.330
&0.309
\\
\hline
\end{array}
}
$
\label{simchi1-3}
\end{table}

\section{Concluding Remarks}

In this paper, we have addressed the problem of the estimation of the mean matrix of the multivariate normal distribution in high dimensional setting. 
Replacing the shrinkage matrix in the Efron-Morris-type estimator with the ridge-type inverse matrix, we have suggested the ridge-type linear shrinkage estimators, where the weighting parameters are estimated by minimizing the Stein unbiased risk estimator under the quadratic loss. 
We have solved the problems of deriving conditions for minimaxity of the suggested linear shrinkage estimators, not only for the known weights, but also for the estimated proposal weights.
The latter problem is clearly harder to solve.
In the framework of the Bayesian model, we have demonstrated that the estimated proposal weights converge to the optimal weights by using random matrix theory.

\medskip
When the constant $c$ of ridge statistic goes to zero, the suggested estimator approaches the Efron-Morris estimator, and  it approaches the James-Stein estimator as $c\to \infty$. 
The choice of $\alh$ greatly affects the estimation accuracy of the proposed estimator. 
The conditions for the minimaxity depend on the choice of $\alh$.
When $n$ is close to $p$ the Efron-Morris estimators are ill-conditioned, but the ridge-type linear shrinkage estimators work well.
This has been supported by the numerical study where we have used the ridge statistic $\alh=\tr(\W)/\min(n-1,p)$.
Although we have considered the two cases $\alh=c$ and $\alh=c\tr(\W)$, it is interesting to obtain the optimal choice of $\alh$, which is a harder problem to be solved in a future.

\medskip
Through out the paper, the covariance matrix $\bSi$ is assumed to be known.
When $\bSi$ is unknown and an estimator $\bSih$ is available, we need to extend the results to the unknown case of $\bSi$.
To this end, similarly to Bodnar et al.(2019), we can consider the bona fide estimators by substituting $\bSih$ into $\bSi$.
Analytical and numerical studies of the bona fide estimators will be done in a future. 

\bigskip
\noindent
{\bf Acknowledgments.}\ \

Research of the first author was supported in part by Grant-in-Aid for JSPS Fellows  (19J22203).
Research of the second author was supported in part by Grant-in-Aid for Scientific Research  (18K11188) from Japan Society for the Promotion of Science.

\renewcommand{\thesection}{\Alph{section}}
\setcounter{section}{0}
\section{Appendix}

\subsection{Unbiased estimation of the risk function}

Consider the case of $n-1\geq p$. 
There exists an $n\times (n-1)$ matrix $\O_1$ such that $\O_1^\top\O_1=\I_{n-1}$ and $\O_1\O_1^\top=\I-n^{-1}\text{\bf 1}_n\text{\bf 1}_n^\top$. 
Let $\Y=\bSi^{-1/2}\X\O_1$ and $\bxi=\bSi^{-1/2}\bTh\O_1$. 
Then,  $\X-\Xb=\X\O_1\O_1^\top$ and $\Y\sim \Nc_{p, n-1}(\bxi, \I_p\otimes \I_{n-1})$.

Let $\H$ and $\L$ be a $p\times p$ orthogonal matrix of eigenvectors and a $p\times p$ diagonal matrix of eigenvalues such that $\Y\Y^\top = \bSi^{-1/2}(\X-\Xb)(\X-\Xb)^\top \bSi^{-1/2} = \H\L\H^\top$, where $\L=\diag (\ell_1, \ldots, \ell_p)$ with the ordering $\ell_1 \geq \cdots \geq \ell_p$.
Let $\bel=(\ell_1, \ldots, \ell_p)^\top$.
Then, we consider the class of the general shrinkage estimators 
\begin{equation}\label{GS}
\bThh_G = \X - \bSi^{1/2}\H\G(\bel)\H^\top\bSi^{-1/2}(\X-\Xb),
\end{equation}
where $\G(\bel)=\diag(g_1(\bel), \ldots, g_p(\bel))$ is a $p\times p$ diagonal matrix of functions of $\bel$.
The estimator $\bThh^{RLS}$ given in (\ref{RLS}) corresponds to $\G(\bel)=a(\L+\alh\I)^{-1}+[b/\tr(\L)]\I$, or 
\begin{equation}\label{gg}
g_i(\bel) = {a \over \ell_i +\alh} + {b\over \tr(\L)}.
\end{equation}

The unbiased estimator of the risk function of $\bThh_G$ under the loss (\ref{loss}) is given from the following theorem. The following way of proof is similar to Stein (1973). The statement is similar to Cand$\rm\grave{e}$s et al (2013) and Hansen (2018).

\begin{thm}
\label{thm:URE}
Let $\Rh(\bThh)$ be an unbiased estimator of the risk function of $\bThh$.
Then, 
$$
\Deh_\G = np\{\Rh(\bThh_G)-\Rh(\X)\}
=\sum_{i=1}^p\ell_i \Big\{g_i^2 - 2(n-1){g_i\over \ell_i}-2\sum_{k\not= i}{g_i-g_k\over \ell_i-\ell_k} - 4{\partial g_i\over \partial \ell_i}\Big\}.
$$
\end{thm}

{\it Proof}.\ \ 
Let $\De = np\{E[L(\bTh, \bThh_G)]-E[L(\bTh, \X)]\}$.
Since $\X-\Xb=\X\P$ for $\P=\I-n^{-1}\text{\bf 1}_n\text{\bf 1}_n^\top$, it can be seen that $\De=E[ \tr\{\H\G^2\H^\top \bSi^{-1/2}\X\P\X^\top\bSi^{-1/2}\}-2\tr\{(\X-\bTh)^\top\bSi^{-1/2}\H\G\H^\top\bSi^{-1/2}\X\P$. 
Let $\Y=\bSi^{-1/2}\X\O_1$ and $\bxi=\bSi^{-1/2}\bTh\O_1$, where $\O_1$ be an $n\times (n-1)$ matrix such that $\O_1^\top\O_1=\I_{n-1}$ and $\O_1\O_1^\top=\I-n^{-1}\text{\bf 1}_n\text{\bf 1}_n^\top$.
Then, $\De$ is rewritten as 
$$
\De=E[\tr(\H\G^2\H^\top \Y\Y^\top)-2\tr\{(\Y-\bxi)^\top\H\G\H^\top\Y\}].
$$
Note that $\Y\sim \Nc_{p, n-1}(\bxi, \I_p\otimes \I_{n-1})$.
Using the Stein identity, we have 
$$
E[\tr\{(\Y-\bxi)^\top\H\G\H^\top\Y\}]= E[\tr\{\bnabla^\top(\H\G\H^\top\Y)\}],
$$ 
where $\bnabla=(\partial/\partial y_{ij})$ for $\Y=(y_{ij})$.
The quantity inside the expectation is calculated as
\begin{align*}
\tr\{\bnabla^\top(\H\G\H^\top\Y)\}
=& \tr[\{\bnabla^\top(\H\G\H^\top)\}\Y] + \tr\{\H\G\H^\top (\bnabla\Y^\top)^\top\} 
\\
=& \tr\{ \Y \bnabla^\top(\H\G\H^\top)\} + (n-1) tr(\G).
\end{align*} 
For $p\times p$ matrix $\U=\U(\Y\Y^\top)$ of functions of $\Y\Y^\top$, it is known that
$$
\bnabla^\top \U = 2 \Y^\top\Dc \U,
$$
where $\Dc=(d_{ij})$ with $d_{ij}=\{(1+\de_{ij})/2\}\partial/\partial s_{ij}$ for $\S=\Y\Y^\top=(s_{ij})$ and the Kronecker's delta $\de_{ij}$.
It is also known that
$$
\Dc(\H\G\H^\top)=\H\G^{(1)}\H^\top,
$$
where $\G^{(1)}=\diag(g_1^{(1)}, \ldots, g_p^{(1)})$ for
$$
g_i^{(1)} = {1\over 2}\sum_{k\not= i} { g_i-g_k \over \ell_i-\ell_k} + {\partial g_i \over \partial \ell_i}.
$$
Thus, one gets
\begin{align*}
\tr\{ \Y \bnabla^\top(\H\G\H^\top)\} =&
2 \tr(\S\H\G^{(1)}\H^\top)=2\tr(\L \G^{(1)})\\
=&
\sum_{i=1}^p\Big\{ \ell_i \sum_{k\not= i} { g_i-g_k \over \ell_i-\ell_k} + 2\ell_i {\partial g_i \over \partial \ell_i}\Big\},
\end{align*}
which gives the expression in Theorem \ref{thm:URE}.
\hfill$\Box$

\subsection{Proof of Theorem \ref{thm:URE0}}

For $\alh$ given in (\ref{alh}), the partial derivative of $\alh$ with respect to $\ell_i$ is
$$
{\partial \alh\over \partial \ell_i} = c_0 = \left\{\begin{array}{ll} 0 & {\rm for}\ \alh=c\\ c &{\rm for}\ \alh=c\tr(\W).\end{array}\right.
$$
Since $g_i=a/(\ell_i+\alh)+b/\{\tr(\L)\}$, we have
$$
{\partial g_i \over \partial \ell_i} = - {a(c_0+1)\over (\ell_i+\alh)^2} - {b \over (\tr\L)^2}.
$$
When $n-1\geq p$, from Theorem \ref{thm:URE}, the unbiased estimator of the risk difference is
\begin{align*}
\Deh =&
\sum_{i=1}^p{\ell_i a^2\over (\ell_i+\alh)^2} + 2 \sum_{i=1}^p{\ell_i ab\over (\ell_i+\alh)\tr\L} + {b^2\over \tr\L}
-2\sum_{i=1}^p{(n-1)a\over \ell_i+\alh}-2{(n-1)pb\over \tr\L}
\\
&
+ 2 \sum_{i=1}^p\sum_{k\not= i}{\ell_i a\over (\ell_i+\alh)(\ell_k+\alh)}
+ 4 \sum_{i=1}^p{(c_0+1)\ell_i a\over (\ell_i+\alh)^2}+4{b\over \tr\L}.
\end{align*}
It is here noted that
\begin{align*}
\sum_{i=1}^p\sum_{k\not= i}{\ell_i \over (\ell_i+\alh)(\ell_k+\alh)}
=&
\sum_{i=1}^p{p-1\over \ell_i+\alh}-\alh \Big(\sum_{i=1}^p {1\over \ell_i+\alh}\Big)^2 + \alh \sum_{i=1}^p {1\over (\ell_i+\alh)^2},
\\
\sum_{i=1}^p{\ell_i \over (\ell_i+\alh)^2} =& \sum_{i=1}^p{1 \over \ell_i+\alh} - \alh\sum_{i=1}^p{1 \over (\ell_i+\alh)^2}.
\end{align*}
Thus, $\Deh$ is rewritten as
\begin{align*}
\Deh =&
\sum_{i=1}^p{\ell_i a^2\over (\ell_i+\alh)^2} + 2 \sum_{i=1}^p{\ell_i ab\over (\ell_i+\alh)\tr\L} + {b^2\over \tr\L}
-2\sum_{i=1}^p{(n-p-1)a\over \ell_i+\alh}-2{(n-1)pb\over \tr\L}
\\
&
-2\alh a \Big(\sum_{i=1}^p {1\over \ell_i+\alh}\Big)^2
+ 2 \sum_{i=1}^p{(2c_0+1)\ell_i a\over (\ell_i+\alh)^2}+4{b\over \tr\L},
\end{align*}
which can be expressed in the form (\ref{URE0}).

In the case of $p>n-1$, it is noted that
\begin{align}
(n-p-1)&\tr((\Y\Y^\top+\alh\I)^{-1})+\alh(\tr((\Y\Y^\top+\alh\I)^{-1}))^2\non\\
=&(p-n+1)\tr((\Y^\top\Y+\alh\I)^{-1})+\alh(\tr((\Y^\top\Y+\alh\I)^{-1}))^2.\label{inveq}
\end{align}
From equation (\ref{inveq}) and Lemma \ref{lem:inv}, the unbiased estimator of the risk difference can be provided by considering the eigenvalue decomposition of $\Y^\top\Y$ and exchanging $p$ and $n-1$ in $\Deh$, and we get the expression in (\ref{URE0-2}).
\hfill$\Box$

\subsection{Useful inequalities}

We here provide some inequalities useful for deriving conditions for minimaxity from the unbiased risk estimators.

\begin{lem}
\label{lem:ineq}
Consider the case of $\alh=c\tr(\W)$ and $n-1\geq p$.
Then, \\

{\rm (1)}\ $1/(1+c)\leq \tr(\V\W)\leq p/(1+cp)$\\

{\rm (2)}\ $p^2/(1+cp)\leq \tr(\V)\tr(\W)\leq p/c$\\

{\rm (3)}\ $\tr(\V^2\W) \leq \tr(\V)/(1+cp) \leq \tr(\V)$\\

{\rm (4)}\ $\tr(\V^3\W) \leq \tr(\V^2)/(1+cp)\leq (\tr\V)^2/(1+cp)$\\

{\rm (5)}\ $\tr(\V^3\W) \leq \tr(\V^2\W)\tr(\V)$.\\

{\rm (6)}\ $\tr(\V^2) \leq (\tr\V)^2[ 1- (p-1)c/\{p(1+c)\}]$.\\

{\rm (7)}\ $\tr(\V^3\W) \leq \tr(\V^2\W)\tr(\V)[ 1- (p-1)c/\{p(1+c)\}]$.\\

{\rm (8)}\ $-c(\tr\V\tr\W-\tr\V\W)/(1+cp) \leq -\tr\V\W+(1+c)\tr\V^2\W^2\leq 0$.

\end{lem}

{\it Proof}.\ \ 
Let $a_i=\ell_i/\tr(\W)$.  
We introduce the expectation $E^*[\cdot]$, defined by $E^*[h(Z)]=p^{-1}\sum_{i=1}^p h(a_i)$ for function $h(\cdot)$.
That is, $Z$ is a random variable with the probability mass $P(Z=a_i)=1/p$ for $i=1, \ldots, p$.
It is noted that $E[Z]=p^{-1}\sum_{i=1}^p a_i=1/p$.

For part (1), $\tr(\V\W)$ is written as $\tr(\V\W)=\sum_{i=1}^p a_i/(a_i+c)=p E[Z/(Z+c)]$.
Since $Z/(Z+c)$ is concave, the Jensen inequality implies that $pE[Z/(Z+c)]\leq p E[Z]/(E[Z]+c)=p/(1+cp)$.
Since $0\leq Z\leq 1$, we have $p E[Z/(Z+c)]\geq p E[Z/(1+c)]=1/(1+c)$.

For part (2), it is noted that $\tr(\V)\tr(\W)=\sum_{i=1}^p 1/(a_i+c)=pE[1/(Z+c)]$.
Since $1/(Z+c)$ is convex, the Jensen inequality gives $pE[1/(Z+c)]\geq p/(E[Z]+c)=p^2/(1+cp)$.
Clearly, $pE[1/(Z+c)]\leq p/c$.

For part (3), $\tr(\V^2\W)=(\tr\W)^{-1}\sum_{i=1}^p a_i/(a_i+c)^2=(\tr\W)^{-1}p E[Z/(Z+c)^2]$.
Since $Z/(Z+c)$ is increasing and $1/(Z+c)$ is decreasing, the covariance inequality implies $E[Z/(Z+c)^2]\leq E[Z/(Z+c)]E[1/(Z+c)]$.
Thus, $\tr(\V^2\W)\leq p^{-1} \tr(\V\W) \tr(\V)$.
From the inequality given in (1), it follows that $\tr(\V^2\W)\leq \tr(\V)/(1+cp)$.

For part (4), the same arguments as in (3) are used to get $\tr(\V^3\W)\leq p^{-1} \tr(\V\W)\tr(\V^2)$. 
From the inequality given in (1), it follows that $\tr(\V^3\W)\leq \tr(\V^2)/(1+cp)$.
Note that $\sum_{i=1}^p b_i^2 \leq (\sum_{i=1}^p b_i)^2$ for nonnegative $b_i$'s.
Hence, $\tr(\V^2)\leq (\tr\V)^2$.

For part (5), it is noted that $\sum_{i=1}^p b_i^3 c_i\leq \sum_{i=1}^p b_i^2 c_i \max(b_1, \ldots, b_p)\leq \sum_{i=1}^p b_i^2 c_i \sum_{i=1}^p b_i$ for nonnegative $b_i$'s and $c_i$'s.
The inequality in (5) follows from this inequality.

For part (6), it is noted that $\sum_{i=1}^p b_i^2=(\sum_{i=1}^p b_i)^2 - \sum\sum_{i\not= j} b_i b_j$.
Then,
\begin{align*}
\tr\V^2=&(\tr\V)^2-\sum\sum_{i\not= j}{1\over ( \ell_i+\alh)( \ell_j+\alh)}
\\
\leq & (\tr\V)^2-\sum\sum_{i\not= j}{1\over ( \ell_i+\alh)( \sum_{k=1}^p \ell_k+\alh)}
\\
=& (\tr\V)^2-{p-1\over 1+c}{\tr(\V)\over \tr(\W)}.
\end{align*}
From the inequality in (2), we have $1/\tr(\W) \geq (c/p)\tr(\V)$, so that $\tr\V^2\leq (\tr\V)^2-(p-1)c/\{p(1+c)\}(\tr\V)^2$.

For part (7), the same arguments as in (6) are used to have
\begin{align*}
\tr\V^3\W=&\tr\V^2\W \tr\V-\sum\sum_{i\not= j}{\ell_i\over ( \ell_i+\alh)^2( \ell_j+\alh)}
\\
\leq & \tr\V^2\W \tr\V-\sum\sum_{i\not= j}{\ell_i\over ( \ell_i+\alh)^2( \sum_{k=1}^p \ell_k+\alh)}
\\
=& \tr\V^2\W \tr\V-{p-1\over 1+c}{\tr(\V^2\W)\over \tr(\W)}.
\end{align*}
From the inequality in (2), we get the requested inequality.

For part (8), $-\tr\V\W+(1+c)\tr\V^2\W^2=-\sum_{i=1}^pa_i/(a_i+c)+(1+c)\sum_{i=1}^pa_i^2/(a_i+c)^2=\sum_{i=1}^pca_i(a_i-1)/(a_i+c)^2=-cpE[Z(1-Z)/(Z+c)^2]\leq0$. Since $Z/(Z+c)$ is increasing and $(1-Z)/(Z+c)$ is decreasing, the covariance inequality implies that $E[Z(1-Z)/(Z+c)^2]\leq E[Z/(Z+c)]E[(1-Z)/(Z+c)]\leq \{E[Z]/(E[Z]+c)\}\{E[1/(Z+c)]-E[Z/(Z+c)]\}=\{1/(1+cp)\}\{E[1/(Z+c)]-E[Z/(Z+c)]\}$. The second inequality holds by Jensen inequality. Thus, $-\tr\V\W+(1+c)\tr\V^2\W^2\geq-c(\tr\V\tr\W-\tr\V\W)/(1+cp)$.
\hfill$\Box$

\subsection{Proof of Theorem \ref{thm:min}}

Consider the case of $g_i=\ah_S/(\ell_i+\alh)$ in Theorem \ref{thm:URE}.
Let $A_0=n-p-1$ for simplicity.
Noting that 
$$
{\partial g_i \over \partial \ell_i} = - {(c_0+1)\ah_S \over (\ell_i+\alh)^2} + {1\over \ell_i+\alh}{\partial \ah_S\over\partial \ell_i},
$$
we have
$$
\sum_{i=1}^p \ell_i {\partial g_i \over \partial\ell_i}= - \sum_{i=1}^p {\ell_i \over (\ell_i+\alh)^2}(c_0+1)\ah_S + \sum_{i=1}^p {\ell_i \over \ell_i+\alh} {\partial \ah_S \over \partial \ell_i}.
$$
Then from Theorem \ref{thm:URE}, the unbiased estimator of the risk difference is written as
\begin{align*}
\Deh_{\ah_S}=&
\sum_{i=1}^p\ell_i\Big({\ah_S\over \ell_i+\alh}\Big)^2 - 2(n-1)\sum_{i=1}^p{\ah_S\over \ell_i+\alh} + 2\ah_S\sum_{i=1}^p\sum_{j\not= i}{\ell_i\over (\ell_i+\alh)(\ell_j+\alh)}\\
&+4(c_0+1)\ah_S\sum_{i=1}^p{\ell_i\over (\ell_i+\alh)^2} - 4\sum_{i=1}^p{\ell_i\over \ell_i+\alh}{\partial\ah_S\over \partial \ell_i}.
\end{align*}
Noting that
\begin{align*}
\sum_{i=1}^p\sum_{j\not= i}{\ell_i\over (\ell_i+\alh)(\ell_j+\alh)}
=& \sum_{i=1}^p{\ell_i\over \ell_i+\alh}\Big\{ \sum_{j=1}^p{1\over \ell_j+\alh} - {1\over \ell_i+\alh}\Big\}\\
=& \sum_{i=1}^p{\ell_i\over \ell_i+\alh}\sum_{j=1}^p{1\over \ell_j+\alh}-\sum_{i=1}^p{\ell_i\over (\ell_i+\alh)^2},
\end{align*}
we get
\begin{align*}
\Deh_{\ah_S}=&
\ah_S^2\tr\V^2\W - 2(n-1)\ah_S\tr\V+2\ah_S\tr(\V\W)\tr\V\\
&-2\ah_S\tr\V^2\W+4(c_0+1)\ah_S\tr\V^2\W  - 4\sum_{i=1}^p{\ell_i\over \ell_i+\alh}{\partial\ah_S\over \partial \ell_i}.
\end{align*}
Note that $\tr\V\W=p-\alh\tr\V$ and $\ah_S=\{A_0\tr\V+\alh(\tr\V)^2\}/\tr\V^2\W-(2c_0+1)$.
Then,
\begin{align*}
- 2(n-1)&\ah_S\tr\V+2\ah_S\tr(\V\W)\tr\V-2\ah_S\tr\V^2\W+4(c_0+1)\ah_S\tr\V^2\W\\
=&- 2\{ (n-1)\ah_S\tr\V-\ah_S(p-\alh\tr\V) \tr\V-2(2c_0+1)\ah_S\tr\V^2\W\}\\
=& -2 \ah_S^2\tr\V^2\W,
\end{align*}
which implies that
\begin{equation}\label{De11}
\Deh_{\ah_S}=
- \ah_S^2\tr\V^2\W  - 4\sum_{i=1}^p{\ell_i\over \ell_i+\alh}{\partial\ah_S\over \partial \ell_i}.
\end{equation}
The partial derivative of $\ah_S$ with respect to $\ell_i$ is 
\begin{align}
{\partial \ah_S\over\partial \ell_i}
=& {1\over \tr\V^2\W}\Big\{- A_0 {1\over(\ell_i+\alh)^2}- A_0\sum_{i=1}^p {c_0\over(\ell_i+\alh)^2}-2\alh{1\over (\ell_i+\alh)^2}\sum_{i=1}^p{1\over \ell_i+\alh}
\non\\
&\hspace{1.8cm}-2\alh\sum_{i=1}^p{c_0\over (\ell_i+\alh)^2}\sum_{i=1}^p{1\over \ell_i+\alh}+c_0\Big(\sum_{i=1}^p{1\over \ell_i+\alh}\Big)^2\Big\}
\non\\
&-{A_0\tr\V+\alh(\tr\V)^2\over (\tr\V^2\W)^2}\Big\{ {1\over (\ell_i+\alh)^2} - 2{\ell_i\over (\ell_i+\alh)^3}- 2\sum_{i=1}^p{c_0\ell_i\over (\ell_i+\alh)^3}\Big\},
\label{ahd}
\end{align}
which yields
\begin{align}
\sum_{i=1}^p& {\ell_i\over \ell_i+\alh}{\partial \ah_S\over \partial \ell_i}
\non\\=&
- {A_0\tr\V^3\W + A_0c_0\tr\V^2\tr\V\W + 2\alh\tr\V\tr\V^3\W + 2c_0\alh\tr\V\tr\V^2\tr\V\W - c_0(\tr\V)^2\tr\V\W \over \tr\V^2\W}\non\\
& - { \{A_0\tr\V+\alh(\tr\V)^2\} \{\tr\V^3\W-2\tr\V^4\W^2-2c_0\tr\V^3\W\tr\V\W\} \over (\tr\V^2\W)^2}
\non\\=&
- {A_0\tr\V^3\W + 2\alh\tr\V\tr\V^3\W + c_0\tr\V\tr\V\W(\alh\tr\V^2-\tr\V) + c_0\tr\V^2\tr\V\W(A_0+\alh\tr\V) \over \tr\V^2\W}\non\\
& - { \{A_0\tr\V+\alh(\tr\V)^2\} \{\tr\V^3\W-2\tr\V^4\W^2-2c_0\tr\V^3\W\tr\V\W\} \over (\tr\V^2\W)^2}.\non
\end{align}
Since $\alh\tr\V^2-\tr\V=-\tr\V^2\W$, from (\ref{De11}), the unbiased estimator of the risk difference is written as
\begin{align}
\Deh_{\ah_S}=& - \tr\V^2\W\Big\{ {A_0\tr\V+\alh(\tr\V)^2\over \tr\V^2\W}-(2c_0+1)\Big\}^2
-4 \sum_{i=1}^p {\ell_i\over \ell_i+\alh}{\partial \ah_S\over \partial \ell_i} \non
\\
=&{\tr\V\over \tr\V^2\W}\Big\{ - (A_0+\alh\tr\V)^2\tr\V + 2(2c_0+1)\tr\V^2\W(A_0+\alh\tr\V) - (2c_0+1)^2{(\tr\V^2\W)^2\over \tr\V}
\non\\
&+ 4A_0{\tr\V^3\W\over \tr\V} + 8\alh\tr\V^3\W - 4 c_0\tr\V\W\tr\V^2\W + 4c_0{\tr\V^2\tr\V\W\over \tr\V}(\A_0+\alh\tr\V)
\non\\
&+4(A_0+\alh\tr\V){\tr\V^3\W-2\tr\V^4\W^2-2c_0\tr\V^3\W\tr\V\W\over \tr\V^2\W}\Big\}. \label{Dehas}
\end{align}
Let $\alh=c$ and $c_0=0$. We have $\tr\V^2\W\leq \tr\V$, $\tr\V^3\W\leq\tr\V^2\leq (\tr\V)^2$ and $\tr\V^3\W/\tr\V^2\W\leq\tr\V$.
For simplicity, let $B=\tr\V$.
Since $\tr\V\W=p-\alh\tr\V$, we obtain the evaluation of $\Deh_{\ah_S}$ as 
\begin{align*}
{\tr\V^2\W\over \tr\V}\Deh_{\ah_S} \leq& B( -A_0^2-2A_0\alh B - \alh^2 B^2 + 2A_0 + 2\alh B 
+ 4A_0+8\alh B + 4A_0 + 4\alh B)
\\
=& (-A_0^2 +10A_0)B + (-2A_0 + 14)\alh B^2 - \alh^2 B^3,
\end{align*}
which is not positive if $A_0=n-p-1\geq 10$.

In the case of $p>n-1$, from \eqref{inveq} and Lemma \ref{lem:inv}, we can get the condition $p-n+1\geq 10$.
These conditions are given in (i) of Theorem \ref{thm:min}.

\medskip
We next obtain condition (ii) of Theorem \ref{thm:min}.
Assume that $n-1\geq p$ and $\alh=c\tr\W$. 
From (\ref{Dehas}), $\Deh_{\ah_S}$ is rewritten as 
\begin{align}
{\tr\V^2\W\over\tr\V}\Deh_{\ah_S}\leq&
(A_0+c\tr\W\tr\V)\tr\V\Big\{ 
 - (A_0+c\tr\W\tr\V) + 2(2c+1){\tr(\V^2\W)\over \tr(\V)}+ 8{\tr\V^3\W\over (\tr\V)^2}
\non\\
&  \hspace{4cm}+ 4c{\tr\V^2\tr\V\W\over (\tr\V)^2}+4{\tr\V^3\W\over \tr\V^2\W\tr\V}\Big\}
\non\\
&
 - (2c+1)^2{(\tr\V^2\W)^2\over \tr\V}-4A_0{\tr\V^3\W\over \tr\V}- 4 c\tr\V\W\tr\V^2\W\non\\
&-8(A_0+c\tr\W\tr\V){\tr\V^4\W^2+c\tr\V^3\W\tr\V\W\over \tr\V^2\W}.
\label{Dehas1}
\end{align}
Using the inequalities in Lemma \ref{lem:ineq}, we can see that 
\begin{align}
 - (A_0+&c\tr\W\tr\V) + 2(2c+1){\tr(\V^2\W)\over \tr(\V)}+ 8{\tr\V^3\W\over (\tr\V)^2}
+ 4c{\tr\V^2\tr\V\W\over (\tr\V)^2}+4{\tr\V^3\W\over \tr\V^2\W\tr\V}
\non\\
\leq &
 - A_0-{cp^2\over 1+cp} + {2(2c+1)\over 1+cp}+ {8\over 1+cp}
+ {4cp\over 1+cp}\Big\{ 1- {(p-1)c\over p(1+c)}\Big\}+4\Big\{ 1- {(p-1)c\over p(1+c)}\Big\}\label{Dehas1sub},
\end{align}
which is not positive if
\begin{equation}
(1+c)(1+cp)A_0 + (p^2-12)c^2 + (p^2-8p-14-4/p)c-14 \geq 0,
\label{scond1}
\end{equation}
or
\begin{equation}
(pA_0+p^2-12)c^2 + \{(p+1)A_0+p^2-8p-14-4/p\}c + A_0-14\geq 0.
\label{scond2}
\end{equation}
Hence from (\ref{Dehas1}) and (\ref{scond1}), the above inequalities are the sufficient conditions given in (ii).

For condition (iii), put $c=1/p$ in (\ref{scond1}).
Then, one gets the sufficient condition
$$
n-p-1 \geq {-p^3+21p^2+14p+16\over 2p(p+1)},
$$
which is given in (iii).
\hfill$\Box$

\subsection{Proof of Theorem \ref{thm:dmin}}

Recall that $g_i=\ah_S/(\ell_i+\alh)+\bh_S/\tr\L$ and $\bh_S = \{(n-1)p-2\}-\tr(\V\W)\ah_S$.
From Theorem \ref{thm:URE} and (\ref{De11}), the unbiased estimator of risk difference can be seen to be
\begin{align*}
\Deh_{\ah_S,\bh_S}
=&-\tr[\V^2\W]\ah_S^2-4\sum{\ell_i\over \ell_i+\alh}{\partial \ah_S\over \partial \ell_i}
+2\ah_S\bh_S\tr\V\W+\bh_S^2-2\{(n-1)p-2\}\bh_S
- {4\over \tr\W}\sum \ell_i{\partial \bh_S\over \partial \ell_i}\\
=&-\tr[\V^2\W]\ah_S^2-4\sum{\ell_i\over \ell_i+\alh}{\partial \ah_S\over \partial \ell_i}-{\bh_S^2\over \tr\W}-{4\over \tr\W}\sum \ell_i{\partial \bh_S\over \partial \ell_i}\\
=&\Deh_{\ah_S,0}-{\bh_S^2\over \tr\W}-{4\over \tr\W}\sum \ell_i{\partial \bh_S\over \partial \ell_i}.
\end{align*}
Note that the partial derivative of $\bh_S$ is
$$
{\partial \bh_S\over \partial \ell_i}
=-{\partial \over \partial \ell_i}\left(\sum_{j=1}^p{\ell_j\over \ell_j+\alh}\ah_S\right)
=-{\ell_i+\alh-\ell_i(1+c_0)\over (\ell_i+\alh)^2}\ah_S-\sum_{j=1}^p{\ell_j\over \ell_j+\alh}{\partial \ah_S\over \partial \ell_i}.
$$
Then it holds that 
\begin{align*}
\sum_{i=1}^p\ell_i{\partial \bh_S\over \partial \ell_i}
=\{-\tr[\V\W]+(1+c_0)\tr[\V^2\W^2]\}\ah_S-\tr[\V\W]\sum_{i=1}^p\ell_i{\partial \ah_S\over \partial \ell_i}.
\end{align*}
From (\ref{ahd}), we can calculate the last term as follows.
\begin{align*}
\sum_{i=1}^p& \ell_i{\partial \ah_S\over \partial \ell_i}
\\=&
- {A_0\tr\V^2\W + A_0c_0\tr\V^2\tr\W + 2\alh\tr\V\tr\V^2\W + 2c_0\alh\tr\V\tr\V^2\tr\W - c_0(\tr\V)^2\tr\W \over \tr\V^2\W}\non\\
& - { \{A_0\tr\V+\alh(\tr\V)^2\} \{\tr\V^2\W-2\tr\V^3\W^2-2c_0\tr\V^3\W\tr\W\} \over (\tr\V^2\W)^2}
\end{align*}
Since $A_0c_0\tr\V^2\tr\W + 2c_0\alh\tr\V\tr\V^2\tr\W - c_0(\tr\V)^2\tr\W=c_0\tr\V\tr\W(\alh\tr\V^2-\tr\V) + c_0\tr\V^2\tr\W(A_0+\alh\tr\V)$, it is rewritten as
\begin{align*}
\sum_{i=1}^p \ell_i{\partial \ah_S\over \partial \ell_i}
=&-\Big\{A_0+2\alh\tr\V+{c_0\tr\V\tr\W(-\tr\V^2\W)\over \tr\V^2\W}+c_0{\tr\V^2\tr\W\over \tr\V^2\W}(A_0+\alh\tr\V)
\\&\quad +(A_0+\alh\tr\V){\tr\V\over \tr\V^2\W}-2(A_0+\alh\tr\V)\tr\V{\tr\V^3\W^2+c_0\tr\W\tr\V^3\W\over(\tr\V^2\W)^2}
\Big\}.
\end{align*}
Since $\alh=c\tr\W$, we have
\begin{align*}
A_0&+2\alh\tr\V-\alh\tr\V+{\alh\tr\V^2\over\tr\V^2\W}(A_0+\alh\tr\V)+(A_0+\alh\tr\V){\tr\V\over\tr\V^2\W}-2(A_0+\alh\tr\V)\tr\V{\tr\V^2\W\over(\tr\V^2\W)^2}
\\=&(A_0+\alh\tr\V)\Big( 1+ {\alh\tr\V^2\over\tr\V^2\W} + {\tr\V\over\tr\V^2\W} -2{\tr\V\over\tr\V^2\W} \Big)\\
=&(A_0+\alh\tr\V)\Big( {\tr\V\over\tr\V^2\W} + {\tr\V\over\tr\V^2\W} -2{\tr\V\over\tr\V^2\W} \Big)=0,
\end{align*}
which shows that $\sum_{i=1}^p \ell_i (\partial \ah_S/\partial \ell_i)=0$.
Hence,  the unbiased estimator of risk difference is expressed as
\begin{equation}
{\tr\V^2\W\over \tr\V}\Deh_{\ah_S,\bh_S}={\tr\V^2\W\over \tr\V}\Deh_{\ah_S,0}-{\tr\V^2\W\over \tr\W \tr\V}\Big[\bh_S^2+4\{-\tr\V\W+(1+c)\tr\V^2\W^2\}\ah_S\Big].
\label{DSD}
\end{equation}
From Lemma \ref{lem:ineq}, $-c(\tr\V\W-\tr\V\W)/(1+cp)\leq-\tr\V\W+(1+c)\tr\V^2\W^2\leq0$, it is evaluated as
\begin{align}
{\tr\V^2\W\over \tr\V}\Deh_{\ah_S,\bh_S}
&\leq {\tr\V^2\W\over \tr\V}\Deh_{\ah_S,0}-{\tr\V^2\W\over \tr\W\tr\V}\Big[4\{-\tr\V\W+(1+c)\tr\V^2\W^2\}\Big\{{\tr\V\over\tr\V^2\W}(A_0+\alh\tr\V)\Big\}\Big]
\non\\
&\leq{\tr\V^2\W\over \tr\V}\Deh_{\ah_S,0}+{4c\over 1+cp}\left(1-{\tr\V\W\over \tr\V\tr\W}\right)\tr\V(A_0+\alh\tr\V)\non\\
&\leq{\tr\V^2\W\over \tr\V}\Deh_{\ah_S,0}+{4c\over 1+cp}\left(1-{1\over p}\right)\tr\V(A_0+\alh\tr\V). \label{Dehasd}
\end{align}

To derive conditions for the minimaxity, from (\ref{Dehas1}) and (\ref{Dehas1sub}), we add the term $4c(1-1/p)/(1+cp)$ to (\ref{Dehas1sub}) in this case. 
Then we can modify conditions (\ref{scond1}) and (\ref{scond2}) by adding $-4(1+c)c(1-1/p)$.
Thus, (\ref{Dehasd}) is negative if
\begin{equation}
(1+c)(1+cp)A_0 + (p^2-16+4/p)c^2 + (p^2-8p-18)c-14 \geq 0,
\label{scond1d}
\end{equation}
or
\begin{equation}
(pA_0+p^2-16+4/p)c^2 + \{(p+1)A_0+p^2-8p-18\}c + A_0-14\geq 0.
\label{scond2d}
\end{equation}
The above inequalities are the sufficient conditions given in (i).

For condition (ii), put $c=1/p$ in (\ref{scond1d}).
Then, one gets the sufficient condition
$$
n-p-1 \geq {-p^4+21p^3+18p^2+16p-4\over 2p^2(p+1)},
$$
which is given in (ii).

\hfill$\Box$

\subsection{Brief description of the random matrix theory}

We use the results of the random matrix theory to prove Theorems \ref{thm:1} and \ref{thm:2} and here we briefly explain the random matrix theory.

Let $(\la_1,\ldots,\la_p)$ be eigenvalues of $\W/(n-1)$. 
The empirical spectral distribution of $\W/(n-1)$ is defined by 
\begin{align*}
F_p(x)=p^{-1}\sum_{i=1}^p1_{[\la_i,\infty)}(x).
\end{align*}
The assumptions (A1)-(A3) in Section \ref{sec:RMT} are sufficient for the existence of the non random limit distribution of $F_p$. 
This property was studied by Mar$\rm{\check{c}}$enko and Pastur (1967), Silverstein and Bai (1995) and Silverstein and Choi (1995), among many others.
Let $F(x)$ be the non random limit distribution of $F_p$, namely
\begin{align*}
F_p(x)\to F(x) \asl p\to \infty.
\end{align*}
For non-decreasing function $G$, the stieltjes transform $m_G$ is given by
\begin{align*}
m_G(z)=\int_{-\infty}^\infty {1\over t-z}dG(t), z\in\C^+\coloneqq\{z\in\C|\Im z>0\}.
\end{align*}
Then, Silverstein and Choi(1995) proved the existence of the following limit of $m_F$, $\check{m}_F(x)\coloneqq\lim_{z\in\C^+\to x}m_F(z)$ for $x\in\R\setminus\{0\}$. They verified some properties of $\check{m}_F(x)$. 
For $x>0$, the following convergence is known as stated in Wang et al.(2015), as $p\to\infty$
\begin{align}
&{1\over p}\tr\left[\left({\W\over n-1}+x\I\right)^{-1}\right]\to \check{m}_F(-x) \as,
\label{Wang1}\\&{1\over p}\tr\left[\left({\W\over n-1}+x\I\right)^{-2}\right]\to \check{m}'_F(-x)={d\check{m}_F(t)\over d t}\Big|_{t=-x} \as. \label{Wang2}
\end{align}

Let $J_p^{(-1)}(z)\coloneqq p^{-1}\tr[(\W/(n-1)-z\I)^{-1}\bPsi^{-1}]$.
Ledoit and P$\rm\acute{e}$ch$\rm\acute{e}$(2011) showed that for $z\in \C^+=\{z\in\C| \Im(z)>0\}$, $\asl p\to\infty$
\begin{align}
J_p^{(-1)}(z)\to J^{(-1)}(z)&=\int_{-\infty}^{\infty}{\tau^{-1} \over \tau(1-c^{-1}-c^{-1}zm_F(z))-z}dH(\tau)\non\\&={m_F(z)\over z}(1-c^{-1}-c^{-1}zm_F(z))-{1\over z}\int_{-\infty}^\infty{dH(\tau)\over \tau}.\label{LP}
\end{align}

Let $z=x+iy$ and $x<0, y>0$, as $y\to 0$, $\tr[(\W/(n-1)-z\I)^{-1}\bPsi^{-1}]\to\tr[(\W/(n-1)-x\I)^{-1}\bPsi^{-1}]$ and $J^{(-1)}(z)=m_F(z)(1-c^{-1}-c^{-1}zm_F(z))/z-(1/z)\int_{-\infty}^\infty1/\tau dH(\tau)\to\check{m}_F(x)(1-c^{-1}-c^{-1}x\check{m}_F(x))/x-(1/x)\int_{-\infty}^\infty1/\tau dH(\tau)$. Then we can extend (\ref{LP}) for $x<0$ as following, $J_p^{(-1)}(x)\coloneqq\tr[(\W/(n-1)-x\I)^{-1}\bPsi^{-1}]\to\check{m}_F(x)(1-c^{-1}-c^{-1}x\check{m}_F(x))/x-(1/x)\int_{-\infty}^\infty1/\tau dH(\tau)\eqqcolon J^{(-1)}(x)$.

By using these convergences, in next subsection, we can show Theorem \ref{thm:1}. 

\subsection{Proof of {\bf Theorem \ref{thm:1}} and Theorem \ref{thm:2}}\label{PfThm1}

For $x<0$, we first show that uniformly convergences of $x\tr[(\W/(n-1)-x\I)^{-1}\Psi^{-1}]/p$, $\tr[(\W/(n-1)-x\I)^{-1}]/p$ and $\tr[(\W/(n-1)-x\I)^{-2}]/p$ in (\ref{Wang1}), (\ref{Wang2}) and (\ref{LP}). The following proof is the same way as Hastie et al.(2019).

Let $f_n(x)=x\tr[(\W/(n-1)-x\I)^{-1}\Psi^{-1}]/p$, it holds that $|f_n(x)|\leq(\min \tau_i)^{-1}\leq c_1^{-1}$. Furthermore $|f'_n(x)|=|\tr[(\W/(n-1)-x\I)^{-2}\W/(n-1)\Psi^{-1}]|/p\leq 8c_2(1+c^{-1/2})^2/c_1^3(1-c^{-1/2})^4\as$ for large enough $n$ by the Bai-Yin theorem. Hence $\{f_n\}$ is equicontinuous. By Arzela-Ascoli theorem, we can verify uniformly convergence of $f_n(x)$.

Let $h_n(x)=\tr[(\W/(n-1)-x\I)^{-1}]/p+x\tr[(\W/(n-1)-x\I)^{-2}]/p=\tr[(\W/(n-1)-x\I)^{-2}\W/(n-1)]$, then $|h_n(x)|\leq 8c_2(1+c^{-1/2})^2/\{c_1^2(1-c^{-1/2})^4\}\as$ for large enough $n$. Furthermore $|h'_n(x)|=2\tr[(\W/(n-1)-x\I)^{-3}\W/(n-1)]/p\leq 32c_2(1+c^{-1/2})^2/\{c_1^3(1-c^{-1/2})^6\}\as$ for large enough $n$. Hence $\{h_n\}$ is equicontinuous. By Arzela-Ascoli theorem, we can verify uniformly convergence of $h_n(x)$.

Let $k_n(x)=(1-p/(n-1))\tr[(\W/(n-1)-x\I)^{-1}]/p-x\{p/(n-1)\}(\tr[(\W/(n-1)-x\I)^{-1}]/p)^2$, it holds that $|k_n(x)|\leq 2/\{c_1(1-c^{-1/2})^2\}\as$ for large enough $n$. Furthermore $|k'_n(x)|\leq 4(1+c+1/c)/\{c_1^2(1-c^{-1/2})^4\} \as$ for large enough $n$. Hence $\{k_n\}$ is equicontinuous. By Arzela-Ascoli theorem, we can verify uniformly convergence of $k_n(x)$.

By these uniformly convergences and continuities of these functions, we have $f_n(-\alh/(n-1))\to -\al J^{(-1)}(-\al)$,  $h_n(-\alh/(n-1))\to \check{m}_F(-\al)-\al\check{m}'_F(-\al)$ and $k_n(-\alh/(n-1))\to (1-1/c)\check{m}_F(-\al)+\al(\check{m}_F(-\al))^2/c \as$ as $n\to \infty$.

By assumption (A2) in Section \ref{sec:RMT}, it is followed that $p^{-1}\tr[\bPsi^{-1}]=p^{-1}\sum\tau_i=\int_{-\infty}^\infty1/\tau dH_p(\tau)\to\int_{-\infty}^\infty 1/\tau dH(\tau)$ as $p\to\infty$. From these convergences, it is seen that as $p\to \infty$

\begin{align*}
{a^*\over n-1}=&{{1\over p}\tr[\bPsi^{-1}]-{\alh/(n-1)\over p}\tr[({\W\over n-1}+{\alh\over n-1}\I)^{-1}\bPsi^{-1}]\over {1\over p}\tr[({\W\over n-1}+{\alh\over n-1}\I)^{-1}]-{\alh/(n-1)\over p}\tr[({\W\over n-1}+{\alh\over n-1}\I)^{-2}]}
\\=&{{1\over p}\tr[\bPsi^{-1}]+f_n(-\alh/(n-1))\over h_n(-\alh/(n-1))}
\\\to&{\int{dH(\tau)\over \tau}-\al J^{(-1)}(-\al)\over\check{m}_F(-\al)-\al\check{m}'_F(-\al)}
\\=&{(1-{1\over c})\check{m}_F(-\al)+{\al\over c}(\check{m}_F(-\al))^2\over\check{m}_F(-\al)-\al\check{m}'_F(-\al)} \as
\end{align*}
and
\begin{align*}
{\ah_{S}\over n-1}&={\left(1-{p\over n-1}\right){1\over p}\tr[({\W\over n-1}+{\alh\over n-1}\I)^{-1}]+{\alh p\over (n-1)^2}({1\over p}\tr[({\W\over n-1}+{\alh\over n-1}\I)^{-1}])^2\over{1\over p}\tr[({\W\over n-1}+{\alh\over n-1}\I)^{-1}]-{\alh/(n-1)\over p}\tr[({\W\over n-1}+{\alh\over n-1}\I)^{-2}]}-{2c_0+1\over n-1}
\\&={k_n(-\alh/(n-1))\over h_n(-\alh/(n-1))}-{2c_0+1\over n-1}
\\&\to{(1-{1\over c})\check{m}_F(-\al)+{\al\over c}(\check{m}_F(-\al))^2\over\check{m}_F(-\al)-\al\check{m}'_F(-\al)} \as
\end{align*}
Therefore, it is concluded that $|\ah_S-a^*|/n\to0 \as$

Since $\tr[\W\Psi^{-1}/(n-1)]/p\to 1 \as \asl n\to \infty$ and $|\tr\W/\{(n-1)p\}-\tr\Psi/p|\to 0 \as \asl n\to\infty$ by assumption (A2) in Section \ref{sec:RMT}, Theorem \ref{thm:2} can be verified.

\end{document}